\documentclass[12pt,a4paper]{article}
\usepackage{amssymb,amsmath}
\usepackage{amsfonts}
\usepackage{amsthm}
\usepackage{a4}
\usepackage{tikz}
\usepackage{enumerate}%
\usepackage{verbatim}
\usepackage{mathalfa}
\usepackage{caption}
\usepackage{float}
\newcommand{\para}{\par\vspace{.25cm}}

\newtheorem*{theorem*}{Theorem}
\newtheorem{theorem}{Theorem}[section]
\newtheorem{lemma}[theorem]{Lemma}

\newtheorem{remark}[theorem]{Remark}
\newtheorem{notation}[theorem]{Notation}

\newcommand{\Q}{\mathbb{Q}}

\newcommand{\C}{\operatorname{Cen}}

\newcommand{\G}{\operatorname{Gal}}

\begin{document}
	\title{\bf Rational group algebras of \\ generalized strongly monomial groups: \\ primitive idempotents and units}
	\author{Gurmeet K. Bakshi$^{a,}${\footnote {Research supported by DST-FIST grant no. SR/FST/MS-II/2019/43 is
	gratefully acknowledged.}}, Jyoti Garg$^{a,}${\footnote {Research supported by Council of Scientific and Industrial Research (CSIR), Govt. of India under the reference no. 09/135(0886)/2019-EMR-I is gratefully acknowledged.}} and Gabriela Olteanu$^{b}$  \\ {\em \small $^a$Department of 
			Mathematics,}\\
		{\em \small Panjab University, Chandigarh 160014, India}\\{\em \small and}\\ {\em \small $^b$Department of Statistics-Forecasts-Mathematics,}\\{\em \small Babe\c s-Bolyai University, Str. T. Mihali 58-60, 400591 Cluj-Napoca, Romania}\\{\em
			\small emails: gkbakshi@pu.ac.in, ~jyotigarg0811@gmail.com, gabriela.olteanu@econ.ubbcluj.ro} }
	\date{}

\maketitle
\begin{abstract} \noindent
We present a method to explicitly compute a complete set of orthogonal primitive idempotents in a simple component with Schur index 1 of a rational group algebra $\mathbb{Q}G$ for $G$ a finite generalized strongly monomial group. For the same groups with no exceptional simple components in  $\mathbb{Q}G$, we describe a subgroup of finite index in the group of units $\mathcal{U}(\mathbb{Z}G)$ of the integral group ring $\mathbb{Z}G$ that is generated by three nilpotent groups for which we give explicit description of their generators. We exemplify the theoretical constructions with a detailed concrete example to illustrate the theory. We also show that the Frobenius groups of odd order with a cyclic complement is a class of generalized strongly monomial groups where the theory developed in this paper is applicable.
\end{abstract}
\para \noindent {\bf MSC2000 :} 16K20, 16S35, 16U60, 20C05, 17C27.
\section{Introduction}
The primitive idempotents in semisimple group algebras have been studied in a series of articles (see \cite{VGO2011, JOdR2012, JOdRVG2013, OVG2015, OVG2016, OVG2022}), where Jespers, Olteanu, del R\'io and Van Gelder described a complete set of orthogonal primitive idempotents in each Wedderburn component 
of a semisimple group algebra $FG$ for various classes of finite groups $G$ and fields $F$, using 
strong Shoda pairs introduced by Olivieri, del R\'io and Sim\'on in \cite{OdRS2004}. This has been done for some classes of groups that include finite nilpotent groups, a class of strongly monomial groups and for fields that are the rationals, finite Galois extensions of the rationals and finite fields. As an application to the knowledge of the primitive idempotents, some of these results were later used to compute
units in integral group rings or minimal non-abelian left group codes (see \cite{JOdR2012,JOdRVG2013,OVG2015}). 

In this article, we compute a complete set of orthogonal primitive idempotents and matrix units in simple components with Schur index 1 of  rational group algebras $\mathbb{Q}G$ for $G$ from a larger class, that is the class of finite generalized strongly monomial groups (see Theorems~\ref{thm3} and \ref{thm4}). In this way, we have extended results to explicitly compute a complete set of primitive idempotents and matrix units in simple components of rational group algebras from classes of groups inside the finite strongly monomial groups to a larger class of groups that is the class of finite generalized strongly monomial groups. The result is proved using different methods to deal with the generalized strongly monomial groups that were introduced and studied in some recent articles (see \cite{BK19}, \cite{BG}, \cite{BGK22}).

	{ In Section 4,  we use the theory developed in Section 3 to study the unit group $\mathcal{U}(\mathbb{Z}G)$ of the integral group ring $\mathbb{Z}G.$ It is well known that for any finite group $G$, $\mathcal{U}(\mathbb{Z}G)$  is finitely generated. The problem of finding a finite generating set of  $\mathcal{U}(\mathbb{Z}G)$ is hard.  To the best of our knowledge, till today it is not known even when $G$ is arbitrary cyclic. However, a lot of progress has been seen in understanding a finite generating set of $\mathcal{U}(\mathbb{Z}G)$ up to a finite index (for the development in this direction, we refer to a recent survey \cite{JS21} and dedicated books \cite{JdR2016}, \cite{Seh93}, \cite{MS} on the topic for history). An early result in this direction is due to Ritter and Sehgal \cite{RS2}, who proved that the Bass cyclic units together with the bicyclic units generate a subgroup of finite index in $\mathcal{U}(\mathbb{Z}G)$ when $G$ is nilpotent of odd order.  In more recent publications (see \cite{JOdR2012, JOdRVG2013}), were given new explicit generators for a subgroup of finite index in $\mathcal{U}(\mathbb{Z}G)$, for $G$ a finite nilpotent group and for $G$ a finite metacyclic group of type $C_{q^m} \rtimes C_{p^n}$, with $C_{p^n}$ acting faithfully on $C_{q^m}$. The advantage of this new method based on the computation of the primitive idempotents is that (i) the proofs were more direct and constructive to obtain explicit sets of generators; (ii) the structure of the subgroup generated by them was given explicitly. In Theorem \ref{cor} of Section 4, we have extended these results to generalized strongly monomial groups $G$ with the property that each simple component in the Wedderburn decomposition of $\Q G$ has Schur index 1 and none of them is exceptional. For such groups $G$, we have explicitly described finitely many generators of three nilpotent groups that together generate a subgroup of finite index in $\mathcal{U}(\mathbb{Z}G)$.}

In Section 5, we have illustrated the theory with detailed and explicit computations of primitive idempotents and units for $G = P \rtimes D_{2^n}$, a semidirect product of the extraspecial $p$-group $P$ of exponent $p$ with the dihedral group of order $2^n$, for $p$ an odd prime, $p \equiv 1 \operatorname{mod} 4$ and $2^{n-1}$  the highest power of $2$ dividing $p-1$. The example for which we have computed explicit sets of primitive orthogonal idempotents and units is complex and of interest, extending previous results regarding this topic. These computations can also be done for groups where the extraspecial $p$-group $P$ has exponent $p^2$ with relative similar computations.

In Section 6, we have shown that the Frobenius groups of odd order with cyclic complement is another class of generalized strongly monomial groups where Theorems \ref{thm3}, \ref{thm4} and \ref{cor} are applicable. Thus for such Frobenius groups, a complete set of primitive idempotents, matrix units of $\Q G$ and the unit group of  $\mathbb{Z}G$ (up to finite index) is computable. It may be pointed out that such Frobenius groups are not known to be strongly monomial.
\section{Preliminaries }\label{prel}
In this section, we introduce some notation and results, mainly from \cite{OdRS2004}, \cite{JdR2016}, \cite{JOdRVG2013}, \cite{BK19} and \cite{BG}.	
Throughout this paper, $G$  denotes a finite group {and all characters of $G$ are complex characters.} A \textit{monomial character} $\chi$ of $G$ is the one which is induced from a linear character of a subgroup of $G$. {In case $\chi$ is a monomial and irreducible character of $G$ which is induced from a linear character $\lambda$ of a subgroup $H$ of $G$,} then the subgroups $H$ and $K= \operatorname{ker} \lambda,$ the kernel of $\lambda,$ satisfy the following:
\begin{itemize}
	\item [(i)] $K \unlhd H$,  $H/K$ is cyclic;
	\item [(ii)] if $g \in G$ and $[H, g] \cap H \subseteq K,$ then $ g \in H.$ Here $[H,g]=\langle g^{-1}h^{-1}gh\, |\,  h \in H \rangle.$\end{itemize}
	{It may be noted that if $H \unlhd N_{G}(K)$, then condition (ii) above is equivalent to saying that $H/K$ is a maximal abelian subgroup of $N_{G}(K)/K$ (see the proof of Proposition 3.5.3 of \cite{JdR2016}).}
\para \noindent A  pair $(H,K)$ of subgroups of $G$ which satisfies (i) and (ii) above is called a \textit{Shoda pair} of $G$ (\cite{OdRS2004}, Definition 1.4) and we call it a Shoda pair of $G$ arising from $\chi$. This shows that  Shoda pairs arise from the monomial irreducible characters of $G$. Conversely, if  $(H,K)$ is a Shoda pair of $G$, by (\cite{CR1962}, Corollary 45.4), we have that $\lambda^G$ is irreducible for any linear character $\lambda$ of $H$ with kernel $K$ and we call it an irreducible character of $G$ arising from the Shoda pair $(H,K).$ 
\para \noindent For an irreducible character $\chi$ of the group $G$,  it is well known {(\cite{TY74}, Proposition 1.1)} that  $${e_{\mathbb{Q}}(\chi) = \frac{\chi(1)}{|G|}\sum_{\sigma \in \operatorname{Gal}(\mathbb{Q}(\chi)/ \mathbb{Q})} {\sum _{g \in G}} \sigma(\chi(g))g^{-1}}$$ {is a  primitive central idempotent of $\Q G$, called the primitive central idempotent of $\Q G$  realized by $\chi$. Here $|G|$ is the order of $G$, $\Q(\chi)$ is the field obtained by adjoining to $\Q$  the character values $\chi(g),~g \in G$,  and $\G(\Q(\chi)/\Q)$ is the Galois group of extension $\Q(\chi)$ over $\Q.$}\para \noindent {Suppose} $(H,K)$ is a Shoda pair of $G$ and $\lambda$, $\lambda'$ are two linear characters of $H$ with kernel $K$. Then clearly $\lambda'= \sigma\circ \lambda$ for some automorphism $\sigma$ of $\Q(\lambda)$ and thus $e_{\mathbb{Q}}(\lambda^G)  = e_{\mathbb{Q}}(\lambda'^G).$  Hence, all the irreducible characters of $G$ which arise from the Shoda pair $(H,K)$ give us the same primitive central idempotent of $\mathbb{Q}G$, which we call the primitive central idempotent of $\mathbb{Q}G$ realized by the Shoda pair $(H,K)$. We call the corresponding simple component $\Q Ge_{\Q}(\lambda^G)$, {{\it the simple component of $\Q G$ realized by the Shoda pair $(H,K)$.}}\para \noindent 
Two Shoda pairs  of $G$ are said to be \textit{equivalent} if they realize the same primitive central idempotent of $\mathbb{Q}G$. The following are some important results on the equivalence of Shoda pairs: \begin{remark}{\label{remark1}} 
	\textit{(i) Two Shoda pairs $(H_{1},K_{1})$, $(H_{2},K_{2})$ of $G$ are equivalent if and only if $H_{1}^g \cap K_{2}= K_{1}^g \cap H_{2}$ for some $g \in G,$ where $H^g=g^{-1}Hg$ (see \cite{JdR2016}, Problem 3.4.3).\\
	(ii) If $H$ is a normal subgroup of $G$ and $(H,K_{1})$, $(H,K_{2})$ are Shoda pairs of $G$, then they are equivalent if and only if $K_{1}$ and $K_{2}$ are conjugate in $G$ (see \cite{JdR2016}, Problem 3.4.4).\\
		(iii) {If $(H_1,K_1)$ and $(H_2,K_2)$ are two Shoda pairs of $G$ and $\lambda_{1}$ (respectively $\lambda_2$) is a linear character of $H_{1}$ (respectively of $H_2$) with kernel $K_{1}$ (respectively $K_2$), then from the definition, it is clear that $(H_1,K_1)$ and $(H_2,K_2)$ are equivalent if and only if  $e_{\Q}(\lambda_{1}^G)= e_{\Q}(\lambda_{2}^G)$ if and only if $\lambda_{1}^G=\sigma \circ \lambda_{2}^G$ for some automorphism $\sigma$ of $\Q(\lambda_{2}^G).$}}
\end{remark}\noindent	
{For a subgroup $H$ of a finite group $G$}, let  $\widehat{H}=\frac{1}{|H|}\displaystyle\sum_{h \in H}h$. Clearly, $\widehat{H}$ is an idempotent of $\mathbb{Q}G$ which is central if and only if $H$ is normal in $G$. For a Shoda pair $(H,K)$ of $G$, let $$\varepsilon(H,K)=\left\{\begin{array}{ll}\widehat{K}, & \hbox{$H=K$;} \\\prod(\widehat{K}-\widehat{L}), & \hbox{otherwise,}\end{array}\right.$$ where $L$ runs over the normal subgroups of $H$ minimal with respect to the property of including $K$ properly. Clearly, $\varepsilon(H,K)$ is an idempotent of the group algebra $\mathbb{Q}G$. Let $e(G,H,K)$ be the sum of all the distinct $G$-conjugates of $\varepsilon(H,K)$, i.e., if $T$ is a right transversal of $\operatorname{Cen}_G(\varepsilon(H,K))$ in $G$, the centralizer of $\varepsilon(H,K)$ in $G$, then $$e(G,H,K)=\sum_{t\in T}\varepsilon(H,K)^t,$$ where $\alpha^g=g^{-1}\alpha g$ for $\alpha\in\mathbb{Q}G$ and $g\in G$. Clearly, $e(G,H,K)$ is a central element of $\mathbb{Q}G$. The element $e(G,H,K)$ is a central idempotent of $\mathbb{Q}G$ if the $G$-conjugates of $\varepsilon(H,K)$ are mutually orthogonal.
\para \noindent
{A Shoda pair $(H,K)$ of $G$ is called a \textit{strong Shoda pair} of $G$ (see \cite{OdRS2004}, Definition 3.1, Proposition 3.3) if the following conditions hold:\begin{description}\item [(a)]$H$ is normal in $\operatorname{Cen}_{G}(\varepsilon(H,K))$;\item [(b)] the distinct $G$-conjugates of $\varepsilon(H,K)$ are mutually orthogonal.\end{description} For reader's convenience, we also recall the following equivalent condition for a pair $(H,K)$ of subgroups of $G$ to be a strong Shoda pair of $G$ (see \cite{JdR2016}, Proposition 3.5.3):
\para \noindent
A pair $(H,K)$ of subgroups of $G$ is a strong Shoda pair of $G$ if and only if the following conditions hold:\\
$(SS1)$ $H \unlhd N_G(K)$,\\
$(SS2)$ $H/K$ is cyclic and a maximal abelian subgroup of $N_G(K)/K$,\\
$(SS3)$ for every $g\in G \backslash N_G(K)$, $\varepsilon(H,K)\varepsilon(H,K)^g = 0.$\para
\noindent If $(H,K)$ is a strong Shoda pair of $G$ and $\lambda$ a linear character of $H$ with kernel $K$, then $e_{\Q}(\lambda^G)=e(G,H,K)$. For the proof, see (\cite{OdRS2004}, Section 3).\para \noindent Recently in \cite{BK19}, the first author and Kaur gave a generalization of the concept of a strong Shoda pair and called them a generalized strong Shoda pair, defined as follows:
\para \noindent If $(H,K)$ is a Shoda pair of $G$ and $\lambda$ is a linear character of $H$ with kernel $K$, then we say that the pair $(H,K)$ is a generalized strong Shoda pair of $G$ if there is a chain $H = H_0 \leq H_1 \leq \cdots \leq H_n = G$ (called a strong inductive chain from $H$ to $G$ of length $n$) of subgroups of $G$ such that the following conditions hold for all $0 \leq i \leq n-1$:
\begin{enumerate}[(i)]
\item $H_i \unlhd \operatorname{Cen}_{H_{i+1}}(e_{\Q}(\lambda^{H_i}))$;
\item the distinct $H_{i+1}$-conjugates of $e_{\Q}(\lambda^{H_i})$ are mutually orthogonal.
\end{enumerate}
 \noindent Note that if $(H,K)$ is a generalized strong Shoda pair with $H=H_0 \leq H_1\leq \cdots \leq H_{n}=G$ as a strong inductive chain, then $(H,K)$ is always a strong Shoda pair of $H_1$. Hence, if a generalized strong Shoda pair $(H,K)$ has a strong inductive chain of length 1, then it is a strong Shoda pair of $G$. Conversely, it is clear that every strong Shoda pair of $G$ is a generalized strong Shoda pair with strong inductive chain $H=H_{0} \leq H_{1}=G$ of length 1.}
{ \begin{remark}\label{remark3}
 (i) If $(H,K) $ is a Shoda pair and $H=H_{0}\leq H_{1}\leq \cdots \leq H_{n}=G$ are subgroups of $G$ such that $H_j \unlhd H_{j+1} $ for some $j$, then the conditions ${\rm (i)}$ and ${\rm (ii)}$ stated above hold true for $i=j$. The condition ${\rm (i)}$ is obvious and ${\rm (ii)}$ holds true because for any $g \in H_{j+1}\backslash \operatorname{Cen}_{H_{j+1}}(e_{\Q}(\lambda^{H_j}))$, $e_{\Q}(\lambda^{H_j})$ and $(e_{\Q}(\lambda^{H_j}))^g$ being distinct primitive central idempotents of $\Q H_{j}$ are mutually orthogonal.\\
 (ii) From part (i), it also follows that if $(H,K)$ is a Shoda pair with $H \unlhd G,$ then $(H,K)$ is a  strong Shoda pair.
 \end{remark}} \noindent A monomial character $\chi$ of $G$ is called {\it generalized strongly monomial} (respectively {\it strongly monomial}) if it arises from a generalized strong Shoda pair of $G$ (respectively strong Shoda pair of $G$). A finite group $G$  is called {\it generalized strongly monomial} (respectively {\it strongly monomial}) if every irreducible character of $G$ is  generalized strongly monomial (respectively strongly monomial).\para \noindent 
{For a ring $R$ with unity, let} $\operatorname{Aut}(R)$  and $\mathcal{U}(R)$ denote the group of automorphisms of $R$ and the group of units of $R$ respectively. Let $\sigma: G \rightarrow \operatorname{Aut}(R)$ and $\tau : G \times G \rightarrow \mathcal{U}(R)$ be two maps which satisfy the following relations,
$$\tau_{gh,x}\sigma_{x}(\tau_{g,h})= \tau_{g,hx} \tau_{h,x}$$ and $$\tau_{g,h}\sigma_{g}(\sigma_{h}(r))= \sigma_{gh}(r) \tau_{g,h},$$for all $g,h,x \in G$ and $r \in R$. Here $\sigma_{g}$ is the image of $g$ under the map $\sigma$ and $\tau_{g,h}$ is the image of $(g,h)$ under the map $\tau.$ {Let $R*_{\tau}^{\sigma} G$ denote the set of finite formal sums $\left\{ \sum z_{g}a_{g} ~|~ a_{g} \in R, g \in G \right\}$, where $z_{g}$ is a symbol corresponding to $g \in G.$ Equality and addition in $R*_{\tau}^{\sigma} G$ are defined componentwise. For $g,h \in G$ and $r \in R,$ by setting
$$z_g z_h = z_{gh} \tau_{g,h},$$ $$rz_g = z_g \sigma_g(r)$$ and extending this rule distributively, $R*_{\tau}^{\sigma}G$ becomes an associative ring, called the {\it crossed product} of $G$ over $R$ with twisting $\tau$ and action $\sigma$. Note that $R*_{\tau}^{\sigma}G$ is a free $R$-module with $\{z_{g}~|~g \in G\}$ an $R$-basis, which we call a basis of units of $R*_{\tau}^{\sigma} G$ as an $R$-module.} A \textit{classical crossed product} is a crossed product $L*_{\tau}^{\sigma} G$, where $L$ is a finite Galois extension of the center  $F$ of $L*_{\tau}^{\sigma}G$, $G=\G(L/F)$ is the Galois group of $L$ over $F$ and $\sigma$ is the natural action of $G$ over $L$. Such a classical crossed product $L*_{\tau}^{\sigma} G$ is commonly denoted by $(L/F,\tau).$\para \noindent  A \textit{cyclic algebra} is a classical crossed product $(L/F,\tau)$, where $L$ over $F$ is a cyclic extension, i.e, $\G(L/F)$ is cyclic. If $\G(L/F)$ is generated by $\phi$ and has order $d$, then $(L/F,\tau)$ is isomorphic to $(L/F,\tau')$, where $\tau'$ is given by $$ \tau'(\phi^i,\phi^j) = \left\{ 
	                                 \begin{array}{ll}
	                                 1, & i+j <d;\\
	                                 z{_\phi}^{d}, &  i+j \geq d,
	                                 \end{array}
\right.
$${where $z_{\phi}\in (L/F,\tau)$ corresponds to $\phi \in \operatorname{Gal}(L/F)$.} Such a cyclic algebra $(L/F,\tau)$ is commonly denoted by $(L/F, \phi, z_{\phi}^d).$ For more details on crossed products, see (Section 8.5 of \cite{MS}, Section 2.6 of \cite{JdR2016}).\para \noindent
\begin{notation}\label{not1} Let $(H,K)$ be a generalized strong Shoda pair of $G$, {$\lambda$ a linear character of $H$ with kernel $K$} and $$H=H_{0}\leq H_{1}\leq \cdots \leq H_{n}=G$$ a strong inductive chain from $H$ to $G$. Throughout the paper, we will use the following notation set in \cite{BK19} and \cite{BG}.
\begin{equation*}
\begin{array}{lll} \noindent
C_{i}& :=& \C_{H_{i+1}}(e_{\Q}(\lambda^{H_{i}})),~ 0 \leq i \leq n-1.\\
k_{i} &:=& \mbox{the index of}~ C_{i}~ \mbox{in}~ H_{i+1} ,~ 0 \leq i \leq n-1.\\
k &:=& k_{0}k_{1} \cdots k_{n-1}.\\ 
\mathtt{k}&:=& \prod_{i=0}^{n-1}|C_{i}/H_{i}|,~\mbox{which is equal to }~\frac{[G:H]}{k}.\\
{(\sigma_{H_{i}})_{x}} &:=& {\mbox{For~} x\in C_{i}/H_{i},~ 0 \leq i \leq n-1, \mbox{~fix}~ \overline{x}\in C_{i}\mbox{~a coset representative of~} x}\\&&
{\mbox{and define the automorphism} ~(\sigma_{H_{i}})_{x} \mbox{~of~}\Q H_{i}e_{\Q}(\lambda^{H_{i}})~ \mbox{by setting}}\\&&
(\sigma_{H_{i}})_{x}(\alpha)=\overline{x}^{-1}\alpha\overline{x},\mbox{~where~} \alpha \in \Q H_{i}e_{\Q}(\lambda^{H_{i}}).\\
{\tau_{H_{i}}} &:=&  {\mbox{~For~}0\leq i\leq n-1,\mbox{the twisting}~\tau_{H_{i}}: C_{i}/H_{i} \times C_{i}/H_{i} \rightarrow \mathcal{U}(\Q H_{i}e_{\Q}(\lambda^{H_{i}}))}~ \\
&& {\mbox{is given by }  \tau_{H_{i}}(x,y) = \overline{xy}^{-1}\overline{x}~\overline{y}e_{\Q}(\lambda^{H_{i}}), \mbox{where}~ x,y \in C_{i}/H_{i}~\mbox{and}~\overline{x},\overline{y},\overline{xy}}\\&&{ \in C_{i}~\mbox{are representatives of the cosets}~x,y,xy~\mbox{respectively}.}\\
\end{array}
\end{equation*}
\noindent For $0\leq i\leq n-1$, we denote by ${\Q H_{i}e_{\Q}(\lambda^{H_{i}})*_{\tau_{H_{i}}}^{\sigma_{H_{i}}}C_{i}/H_{i}}$ the crossed product of $\Q H_{i}e_{\Q}(\lambda^{H_{i}})$ over $C_{i}/H_{i}$ with action $\sigma_{H_{i}}:C_{i}/H_{i}\rightarrow \operatorname{Aut}(\Q H_{i}e_{\Q}(\lambda^{H_{i}}))$, where $\sigma_{H_{i}}(x)=(\sigma_{H_{i}})_{x},~x \in C_{i}/H_{i}$ and the twisting $\tau_{H_{i}}: C_{i}/H_{i} \times C_{i}/H_{i} \rightarrow \mathcal{U}(\Q H_{i}e_{\Q}(\lambda^{H_{i}}))$ defined as above.
\end{notation}
\noindent 
{It is proved in (\cite{BK19}, Proof of  Theorem 3) that for all $0 \leq i \leq n-1,$ \begin{equation}{\label{eq1}} \theta_i:
  			\Q H_{i+1}e_{\mathbb{Q}}(\lambda^{H_{i+1}}) \to  M_{k_i}(\Q C_{i}e_{\mathbb{Q}}(\lambda^{H_{i}}) )
  		\end{equation}  given by $\alpha \mapsto (\alpha_{ml})_{k_{i} \times k_{i}}$, where  $\alpha \in \Q H_{i+1}e_{\Q}(\lambda^{H_{i+1}})$,  $\alpha_{ml}= e_{\Q}(\lambda^{H_{i}})t_{(i,l)} \alpha t_{(i,m)}^{-1}e_{\Q}(\lambda^{H_{i}})$ is an isomorphism.  However,  $\Q C_{i}e_{\mathbb{Q}}(\lambda^{H_{i}})$ is naturally identified with the crossed product  $\Q H_{i}e_{\Q}(\lambda^{H_{i}})*_{\tau_{H_{i}}}^{\sigma_{H_{i}}}C_{i}/H_{i}$,  and one obtains the following isomorphism  for $0\leq i \leq n-1,$ \begin{equation}\label{eq -1}
  			\Q H_{i+1}e_{\Q}(\lambda^{H_{i+1}})\stackrel{\theta_{i}}{\cong} M_{k_{i}}(\Q H_{i}e_{\Q}(\lambda^{H_{i}})*_{\tau_{H_{i}}}^{\sigma_{H_{i}}}C_{i}/H_{i}). \end{equation} A recursive application of the  isomorphisms $\theta_i, \theta_{i-1}, \cdots, \theta_0$ yield the following isomorphism for all $0\leq i \leq n-1$: \begin{equation}{\label{eq 0}}
  			{\footnotesize \mathbb{Q}H_{i+1}e_{\mathbb{Q}}(\lambda^{H_{i+1}}) \cong  {M_{k_{i}}(\cdots(M_{k_{0}}(\Q H e_{\Q}(\lambda)
  					*^{\sigma_{H_{0}}}_{\tau_{H_{0}}} C_{0}/H_{0}) *^{\sigma_{H_{1}}}_{\tau_{H_{1}}} \cdots) *^{\sigma_{H_{i}}}_{\tau_{H_{i}}} C_{i}/H_{i})}.} \end{equation} Let us denote the above isomorphism by ${\Theta}_i$. The particular case when $i=n-1$  yields the  following result on the structure of a simple component of $\mathbb{Q}G$ realized by a generalized strong Shoda pair :}\para
\begin{theorem}\label{thm1}{\rm  (\cite{BK19}, Theorem 3)}  Let $(H,K)$ be a generalized strong Shoda pair of $G$ and {$\lambda$ a linear character of $H$ with kernel $K$.} Let $H=H_{0}\leq H_{1}\leq \cdots \leq H_{n}=G$ be a strong inductive chain from $H$ to $G$. Then {\footnotesize $$\mathbb{Q}Ge_{\mathbb{Q}}(\lambda^{G}) \cong  {M_{k_{n-1}}(\cdots(M_{k_{0}}(\Q H e_{\Q}(\lambda)
			*^{\sigma_{H_{0}}}_{\tau_{H_{0}}} C_{0}/H_{0}) *^{\sigma_{H_{1}}}_{\tau_{H_{1}}} \cdots) *^{\sigma_{H_{n-1}}}_{\tau_{H_{n-1}}} C_{n-1}/H_{n-1})},$$}where $C_{i},~\sigma_{H_{i}},~\tau_{H_{i}},~k_{i}$ are as defined above. \end{theorem}
\noindent  In \cite{BG}, the first two authors explained the structure of the algebra\linebreak $ M_{k_{n-1}}(\cdots(M_{k_{0}}(\Q He_{\mathbb{Q}}(\lambda) *^{\sigma_{H_{0}}}_{\tau_{H_{0}}} C_{0}/H_{0}) *^{\sigma_{H_{1}}}_{\tau_{H_{1}}} \cdots)*^{\sigma_{H_{n-1}}}_{\tau_{H_{n-1}}} C_{n-1}/H_{n-1})$. We will briefly recall it below. 
\begin{notation}\label{not2} {For a generalized strong Shoda pair $(H,K)$ of $G$, $\lambda$ a linear character of $H$ with kernel $K$ and $$H=H_{0}\leq H_{1}\leq \cdots \leq H_{n}=G$$ a strong inductive chain from $H$ to $G$, the following notation introduced in \cite{BG} will also be used throughout.}
\begin{equation*}
\begin{array}{lll} \noindent
\mathcal{A} &:=& M_{k_{n-1}}(\cdots(M_{k_{0}}(\Q He_{\Q}(\lambda) *^{\sigma_{H_{0}}}_{\tau_{H_{0}}} C_{0}/H_{0}) *^{\sigma_{H_{1}}}_{\tau_{H_{1}}} \cdots)*^{\sigma_{H_{n-1}}}_{\tau_{H_{n-1}}} C_{n-1}/H_{n-1}).\\
\mathbb{F} &:=& \mathcal{Z}(\mathcal{A}), \mbox{the center of~} \mathcal{A}.\\
\mbox{E} &:=& \Q He_{\Q}(\lambda).\\
\mathbb{E} &:=& k \times k~ \mbox{scalar matrices}~ \operatorname{diag}(\alpha,\alpha,\cdots,\alpha)_{k}, \alpha \in \mbox{E}.\\
\mbox{F} &:=&\{\alpha \in \mbox{E} ~|~ \operatorname{diag}(\alpha,\alpha,\cdots,\alpha)_{k} \in \mathbb{F}\}.\\
\mathcal{B} &:=& M_{{k}}(\mbox{F}).\\
\operatorname{Cen}_{\mathcal{A}}(\mathcal{B})&:=& \mbox{the centralizer of~} \mathcal{B} \mbox{~in~} \mathcal{A}.\\
 \mathcal{G} &:=& \operatorname{Gal}(\mathbb{E}/\mathbb{F}) .\\
E_{ij} &:=& \mbox{the matrix in}~\mathcal{A}~\mbox{whose}~ $(i,j)$\mbox{-th entry is}~\varepsilon(H,K)~\mbox{and all other entries}\\&& \mbox{are zero,}~{1\leq i,j\leq k.}\\
\zeta_{[H:K]} &:=& \mbox{a primitive}~[H:K]\mbox{-th root of unity.} 
\end{array}
\end{equation*}
\end{notation}\noindent
{A group $G$  is said to be of type $G_{1}$-by-$G_{2}$ if $G$ has a normal subgroup isomorphic to $G_{1}$ and the quotient by $G_{1}$ is isomorphic to $G_{2}$.}
{\begin{theorem}{\label{thm2}}{\rm (\cite{BG}, Section 3)}
	Let $(H,K)$ be a generalized strong Shoda pair of $G$, $\lambda$ a linear character of $H$ with kernel $K$ and $H=H_{0}\leq H_{1} \leq \cdots \leq H_{n}=G$ a strong inductive chain from $H$ to $G$. Then 
	\begin{enumerate}[(i)]
	\item $\Q Ge_{\Q}(\lambda^G) \cong \mathcal{B} \otimes_{\mathbb{F}} \operatorname{Cen}_{\mathcal{A}}(\mathcal{B}) \cong M_{k}(\operatorname{Cen}_{\mathcal{A}}(\mathcal{B})).$
	\item $\mathbb{E}$ is a Galois extension of $\mathbb{F}$, $\operatorname{dim}_{\mathbb{F}}(\mathbb{E})=\mathtt{k}=\prod_{i=0}^{n-1}|C_{i}/H_{i}|$, $\operatorname{Cen}_{\mathcal{A}}(\mathcal{B})$ contains $\mathbb{E}$ as a maximal subfield and $\operatorname{dim}_{\mathbb{F}}(\operatorname{Cen}_{\mathcal{A}}(\mathcal{B}))=(\operatorname{dim}_{\mathbb{F}}(\mathbb{E}))^2=\mathtt{k}^2.$
	\item there exist units $\{z_{\sigma}~|~\sigma \in \mathcal{G}\}$ in $\operatorname{Cen}_{\mathcal{A}}(\mathcal{B})$, which form a basis of $\operatorname{Cen}_{\mathcal{A}}(\mathcal{B})$ as a vector space over $\mathbb{E}$ and $\operatorname{Cen}_{\mathcal{A}}(\mathcal{B})\cong (\mathbb{E}/\mathbb{F},\tau),$ where $\tau: \mathcal{G}\times \mathcal{G} \rightarrow \mathcal{U}(\mathbb{E})$ is given by
	\begin{equation}{\label{n}}
	\tau(\sigma,\sigma')=z_{\sigma\sigma'}^{-1}z_{\sigma}z_{\sigma'}.
	\end{equation} 
	\item $\mathcal{G}$ is of the type $(((C_{0}/H_{0}$-by-$C_{1}/H_{1})$-by-$C_{2}/H_{2})\cdots)$-by-$C_{n-1}/H_{n-1}.$ In particular, the order of $\mathcal{G}$ is $\mathtt{k}$. 
	\end{enumerate}
\end{theorem}}
\begin{remark}{\label{rr}}
(a) If $(H,K)$ is a strong Shoda pair of $G$, then $E=\Q H\varepsilon(H,K)$ and $F=\Q H\varepsilon(H,K)^{N_{G}(K)/H}$, the fixed subfield of $\Q H\varepsilon(H,K)$ under the action of $N_{G}(K)/H$ by conjugation. Also, $\mathcal{B}\cong M_{k_{0}}(\Q H\varepsilon(H,K)^{N_{G}(K)/H})$ and $\operatorname{Cen}_{\mathcal{A}}(\mathcal{B})\cong \Q H\varepsilon(H,K)*_{\tau}^{\sigma}N_{G}(K)/H$. Further any set of coset representatives of $N_{G}(K)/H$ can be taken as a basis of units of $\operatorname{Cen}_{\mathcal{A}}(\mathcal{B})$ as $E$-vector space.\\
(b) If $(H,K)$ is a generalized strong Shoda pair of $G$ which is not a strong Shoda pair, the task to find a {basis of units} {$\{z_{\sigma}~|~\sigma\in\mathcal{G}\}$ of $\operatorname{Cen}_{\mathcal{A}}(\mathcal{B})$ as $\mathbb{E}$-vector space is a bit more complicated}, it is explained in Section 3.1 of \cite{BG} and illustrated with examples in the same paper.
\end{remark}    \noindent		
If $E/F$ is a finite Galois extension, then there exists an element $w\in E$ such that $\{\sigma(w)~|~\sigma \in \G(E/F)\}$ is an $F$-basis of $E$. Such a basis is called {\it a normal basis} and $w$ {\it a normal element of} $E/F.$\para \noindent	{Let $R$ be a semisimple ring, i.e., R is a direct sum of a finite number of minimal left ideals.  Recall that a {\it complete set of orthogonal primitive
idempotents} of $R$ is  a family $\{e_1, \ldots , e_t \}$ of elements of $R$  such that: (i) each $e_i$ is a non-zero idempotent element; (ii) if $i\neq  j$, then $e_i e_j = 0$; (iii) $e_1+\cdots+e_t=1$; (iv) none of the $e_i,~1\leq i\leq t,$ can be written as $e_i = e_i'+e_i''$ , where $e_i'$, $e_i''$ are non zero idempotents with $e_i' e_i''=0$.}	\para \noindent
{It may also be useful to recall that {\it a complete set of matrix units} of a matrix algebra $M_{n}(D)$ over a division ring $D$ is a set of elements $E_{ij}\in M_{n}(D)$, with $1\leq i,j \leq n$ such that  $E_{ij}E_{kl}=\delta_{jk}E_{il}$ and $E_{11}+\cdots+E_{nn}=I_{n}$, where $\delta_{ij}$ is the Kronecker delta function and 
$I_{n}$ is the identity matrix. Note that $E_{ii},~1\leq i\leq n,$ form a complete set of orthogonal primitive idempotents of the simple ring $M_{n}(D)$.}
\para \noindent In \cite{JOdRVG2013}, there are three results regarding the computation of complete sets of orthogonal primitive idempotents (Theorem 4.1), matrix units (Corollary 4.2) and description of 
a subgroup of finite index in the unit group of the integral group ring (Theorem 5.4) for certain classes of strongly monomial groups. In the next two sections, the aim is to extend these results to generalized strongly monomial groups. 
\section{Primitive Idempotents }
In this section, we will explicitly construct a complete set of orthogonal primitive idempotents of $\Q Ge_{\Q}(\lambda^G)$ in the particular case when this simple component has Schur index  1. Recall that the {\it Schur index} of a central simple algebra $A$ is defined as the degree of the division algebra $D$ once we write $A$ as isomorphic to some $\mathcal{M}_m(D)$. By the degree of the division algebra $D$, we mean the square root of the dimension of $D$ over its center.  \para \noindent {From Theorems \ref{thm1} and \ref{thm2}, it is clear that in order to find a complete set of orthogonal primitive idempotents of $\Q Ge_{\Q}(\lambda^G)$, it is enough to find primitive idempotents of the central simple $\mathbb{F}$-subalgebra of $\Q Ge_{\Q}(\lambda^G)$, which maps to $\mathcal{B}$ and $\operatorname{Cen}_{\mathcal{A}}(\mathcal{B})$, under the isomorphism $\Q Ge_{\Q}(\lambda^G) \cong \mathcal{B} \otimes_{\mathbb{F}} \operatorname{Cen}_{\mathcal{A}}(\mathcal{B})$. For the convenience of language, we will denote the central simple $\mathbb{F}$-subalgebra of  $\Q Ge_{\Q}(\lambda^G)$ which maps to $\mathcal{B}$ (respectively $\operatorname{Cen}_{\mathcal{A}}(\mathcal{B})$) again by $\mathcal{B}$ (respectively $\operatorname{Cen}_{\mathcal{A}}(\mathcal{B})$).}
The first step is to  find the primitive idempotents of the $\mathbb{F}$-subalgebra $\mathcal{B}$ of  $\Q Ge_{\Q}(\lambda^G)$ for which no assumption on the Schur index is required. We use Notations~\ref{not1} and \ref{not2}.
\begin{lemma}{\label{l1}}
  Let $(H,K)$ be a generalized strong Shoda pair, {$\lambda$ a linear character of $H$ with kernel $K$ and} $H=H_{0}\leq H_{1} \leq \cdots \leq H_{n}=G$ a strong inductive chain from $H$ to $G$. For $0 \leq i \leq n-1$, let $T_{i}=\{t_{(i,1)},t_{(i,2)}, \ldots, t_{(i,k_{i})} \}$ be a left transversal of $C_{i}$ in $H_{i+1}$ {with $t_{(i,1)}=1$}. Let $T= T_{0}T_{1}\cdots T_{n-1}=\{t_{(0,s_{0})}t_{(1,s_{1})}\cdots t_{(n-1,s_{n-1})}~|~t_{(i,s_{i})} \in T_{i},~0 \leq i \leq n-1,~1 \leq s_{i} \leq k_{i} \}$. Then $$\{ t^{-1} \,\varepsilon(H,K)\,t~|~ t \in T\}$$ is a complete set of orthogonal primitive idempotents of $\mathcal{B}$. 
\end{lemma}
{\noindent {\bf Proof.} For $0\leq i \leq n-1$, recall from equations (\ref{eq1}) and (\ref{eq -1}), the isomorphism: $$ \theta_i:
  	\Q H_{i+1}e_{\mathbb{Q}}(\lambda^{H_{i+1}}) \to  M_{k_i}(\Q H_{i}e_{\mathbb{Q}}(\lambda^{H_{i}}) *^{\sigma _{H_{i}}}_{\tau_{H_{i}}} C_i/H_i),
  	$$  given by $\alpha \mapsto (\alpha_{ml})_{k_{i}}$, where  $\alpha \in \Q H_{i+1}e_{\Q}(\lambda^{H_{i+1}})$,  $\alpha_{ml}= e_{\Q}(\lambda^{H_{i}})t_{(i,l)} \alpha t_{(i,m)}^{-1}e_{\Q}(\lambda^{H_{i}})$.  Also from equation (\ref{eq 0}), we have the isomorphism  
  	$$\Theta_i: \Q H_{i+1}e_{\Q}(\lambda^{H_{i+1}}) \to {M_{k_{i}}(\ldots(M_{k_{0}}(\Q H e_{\Q}(\lambda)
  		*^{\sigma_{H_{0}}}_{\tau_{H_{0}}} C_{0}/H_{0}) *^{\sigma_{H_{1}}}_{\tau_{H_{1}}} \ldots) *^{\sigma_{H_{i}}}_{\tau_{H_{i}}} C_{i}/H_{i})}.$$ 
  \noindent  Consider  $t=t_{(0,s_{0})} t_{(1,s_{1})}\cdots t_{(n-1,s_{n-1})} \in T$, where $t_{(i,s_{i})} \in T_{i},$  $ 1 \leq s_i\leq k_i$. Let $\alpha =t^{-1}\varepsilon(H,K) t$.  For $0\leq i \leq n-1$, let  us denote $ t_{(0,s_{0})} t_{(1,s_{1})}\cdots t_{(i,s_{i})}$ by $t_i $, and  $e_{\Q}(\lambda^{H_{i+1}})t_i^{-1}\varepsilon(H,K) t_i e_{\Q}(\lambda^{H_{i+1}})$ by $\alpha_i$, so that $t_{n-1}=t$ and $\alpha_{n-1}= \alpha$. Observe that $\alpha_i \in \Q H_{i+1}e_{\Q}(\lambda^{H_{i+1}}) $ for all $i$. \vspace{.2cm} \\
  	{\bf Step 1:} For any $0 \leq i \leq n-1$, we show that $\theta_{i}(\alpha_i) $ is a $k_i \times k_i$ matrix with $\alpha_{i-1}$ at the $(s_{i}, s_{i})$ entry  and all other entries zero. By $\alpha_{-1}$, we mean $\varepsilon(H,K)$. \vspace{.2cm} \\
  	From the definition $$\theta_{i}(\alpha_i) = (\alpha_{ml}),$$ where  \begin{equation}
  		\begin{array}{lll}\noindent \alpha_{ml}&=&   e_{\Q}(\lambda^{H_{i}})t_{(i, l)}e_{\Q}(\lambda^{H_{i+1}})t_i^{-1}\varepsilon(H, K) t_ie_{\Q}(\lambda^{H_{i+1}}) t_{(i,m)}^{-1}e_{\Q}(\lambda^{H_{i}})\\ &=& e_{\Q}(\lambda^{H_{i}})t_{(i,l)}t_i^{-1}\varepsilon (H,K) t_i t_{(i,m)}^{-1}e_{\Q}(\lambda^{H_{i}})\\ &=& e_{\Q}(\lambda^{H_{i}})t_{(i,l)}t_{(i,s_i)}^{-1}t_{i-1}^{-1}\varepsilon(H,K) t_{i-1} t_{(i,s_i)}t_{(i,m)}^{-1}e_{\Q}(\lambda^{H_{i}}), ~{\rm as}~ t_i= t_{i-1}t_{(i,s_i)}\\&=& e_{\Q}(\lambda^{H_{i}})t_{(i,l)}t_{(i,s_i)}^{-1}e_{\Q}(\lambda^{H_{i}})t_{i-1}^{-1}\varepsilon(H,K) t_{i-1} e_{\Q}(\lambda^{H_{i}})t_{(i,s_i)}t_{(i,m)}^{-1}e_{\Q}(\lambda^{H_{i}}).   \end{array} 
\end{equation} 
  		 The second equality follows because $e_{\Q}(\lambda^{H_{i+1}})$ commutes with $t_{(i,l)}$, $ t_{(i,m)}$ and from (\cite{BK19}, Lemma 1) $e_{\Q}(\lambda^{H_{i}})e_{\Q}(\lambda^{H_{i+1}})= e_{\Q}(\lambda^{H_{i}})= e_{\Q}(\lambda^{H_{i+1}})e_{\Q}(\lambda^{H_{i}})$. The last equality is because  $ \varepsilon(H,K) = \varepsilon(H,K)e_{\Q}(\lambda^{H_{i}})= e_{\Q}(\lambda^{H_{i}}) \varepsilon(H,K)$  and $ e_{\Q}(\lambda^{H_{i}})$ commutes with $t_{i-1}$. If $ l \neq s_i$, then  $e_{\Q}(\lambda^{H_{i}})t_{(i,l)}t_{(i,s_i)}^{-1}e_{\Q}(\lambda^{H_{i}})t_{(i,s_i)}t_{(i,l)}^{-1} =0$ as  the distinct $H_{i+1}$- conjugates of  $e_{\Q}(\lambda^{H_{i}})$ are mutually orthogonal. Consequently $e_{\Q}(\lambda^{H_{i}})t_{(i,l)}t_{(i,s_i)}^{-1}e_{\Q}(\lambda^{H_{i}})=0$, and hence $\alpha_{ml}=0$.   Similarly if $ m \neq s_i$, then $e_{\Q}(\lambda^{H_{i}})t_{(i, s_i)}t_{(i,m)}^{-1}e_{\Q}(\lambda^{H_{i}})$ is zero and hence $\alpha_{ml}=0$.  If $l=m = s_i$,  it is clear that then $\alpha_{ml} = \alpha_{i-1}$.  This proves step 1. 
  	\noindent \vspace{.2cm} \\{\bf Step 2:} We now show that   $\alpha = t^{-1}\varepsilon(H,K) t$ belongs to $\mathcal{B}$ and it is a primitive idempotent. 
\noindent \vspace{.2cm} \\For this step, we are actually required to show that $\Theta_{n-1}(\alpha) \in M_k(F)$ and it is a primitive idempotent.  For any $m \geq 1$, denote by $E^{(m)}_{ii}$ the $m\times m$ matrix whose $(i,i)$ entry is $\varepsilon(H,K)$ and all other entries are zero. From step 1, it follows that, for any $ 0\leq i \leq n-1$, $\Theta_i(\alpha_i)$ is a $k_i\times k_i$ block matrix with each block of size $k_0k_1\cdots k_{i-1} \times k_0k_1\cdots k_{i-1} $ such that the $(s_i, s_i)$ block is $\Theta_{i-1}(\alpha_{i-1})$ are all other blocks are zero.  Recall that $\alpha_{-1}=\varepsilon(H,K)$. By taking $i=0,1, \cdots, n-1$, it yields that  $\Theta_{0}(\alpha_0)= E^{(k_o)}_{ii}$  with $i=s_0$, $\Theta_{1}(\alpha_1)= E^{(k_ok_1)}_{ii}$  with $i=(s_1-1)k_0+s_0$,  $\Theta_{2}(\alpha_2)= E^{(k_ok_1k_2)}_{ii}$  with $i=(s_2-1)k_0k_1+(s_1-1)k_0+s_0$  and finally,  $\Theta_{n-1}(\alpha_{n-1})= E^{(k_ok_1\cdots k_{n-1})}_{ii}$  with  $i=\sum_{j=1}^{n-1} (s_{j}-1)k_{0}\cdots k_{j-1}+s_{0}$.  
  	As $ k_0k_1\cdots k_{n-1}=k$, $ \alpha_{n-1}= \alpha$ and $\varepsilon(H,K)$ is the identity of $F$,  the desired step  is proved. 
  \noindent \vspace{.2cm} \\{\bf Step 3:} We will see that if $t, t' \in T$, $t \neq t'$, then  $t^{-1}\varepsilon(H,K) t$ and $t'^{-1}\varepsilon(H,K) t'$ are orthogonal. \vspace{.2cm} \\ 
  Write $ t=t_{(0,s_{0})} \cdots t_{(n-1,s_{n-1})}$ and $t'=t_{(0,s'_{0})} \cdots t_{(n-1,s'_{n-1})}$, where $1 \leq s_{i},s'_{i} \leq k_{i}$. If $t \neq t'$, then choose the least integer $m$ such that $s_{m}\neq s_{m}'$. In case $t^{-1}\varepsilon(H,K)t$ and $t'^{-1}\varepsilon(H,K)t'$ are not orthogonal, then they are equal and hence $$\sum_{j=1}^{n-1} (s_{j}-1)k_{0}\cdots k_{j-1}+s_{0} = \sum_{j=1}^{n-1} (s'_{j}-1)k_{0}\cdots k_{j-1}+s'_{0}.$$ Since $s_i=s_i'$ for $0\leq i\leq m-1$, we obtain that $\sum_{j=m+1}^{n-1} (s_{j}-1)k_{0}\cdots k_{j-1}+(s_{m}-1)k_{0}\cdots k_{m-1} = \sum_{j=m+1}^{n-1} (s'_{j}-1)k_{0}\cdots k_{j-1}+(s'_{m}-1)k_{0}\cdots k_{m-1}$, which gives $\sum_{j=m+1}^{n-1} (s_{j}-s'_{j})k_{0}\cdots k_{j-1}=(s'_{m}-s_{m})k_{0}\cdots k_{m-1}$. Now, $k_{0}\cdots k_{m}$ divides left hand side of the previous equation, so it also divides right hand side. This implies $k_{m}| s_{m}-s'_{m}$, which is a contradiction as $|s_{m}-s'_{m}|<k_{m}$. 
  \noindent \vspace{.2cm} \\{\bf Step 4:} Finally, $\{ t^{-1} \,\varepsilon(H,K)\,t~|~ t \in T\}$ is a complete set of orthogonal primitive idempotents of $\mathcal{B}.$  \vspace{.2cm} \\ This follows immediately from step 3 because there are precisely $k$ primitive idempotents of  $M_k(F)$ and the cardinality of $T$ is $k_ok_1\cdots k_{n-1}=k$.  The lemma is  thus proved. \qed}
\para \noindent {Our next  step is to find the primitive idempotents of $\operatorname{Cen}_{\mathcal{A}}(\mathcal{B})$  for which we shall require an assumption on the Schur index. Before we begin with its statement, for reader's convenience, we briefly recall from \cite{BG}, the existence of the units $\{z_{\sigma}~|~\sigma \in \mathcal{G}\}$ in $\operatorname{Cen}_{\mathcal{A}}(\mathcal{B})$, which form a basis of $\operatorname{Cen}_{\mathcal{A}}(\mathcal{B})$ as a vector space over $\mathbb{E}$. In (\cite{BG}, Proof of step (iv), Page 9), we have mentioned that, by Noether-Skolem Theorem (\cite{JdR2016}, Theorem 2.1.9), for every $\sigma \in \G(\mathbb{E}/\mathbb{F}),$ there exists an invertible $z_{\sigma}$ in $\C_{\mathcal{A}}(\mathcal{B})$ such that $\sigma(\alpha)= z_{\sigma}^{-1} \alpha z_{\sigma}$ for all $\alpha \in \mathbb{E}.$ We can quickly see that $\{z_{\sigma} ~|~ \sigma \in \G(\mathbb{E}/\mathbb{F})\}$ are linearly independent over $\mathbb{E}$. Suppose $k$ is the least integer with $\sigma_{1},\sigma_{2},\ldots,\sigma_{k}$ in $\mathcal{G}$ such that $\sum_{i=1}^{k}a_{\sigma_{i}}z_{\sigma_{i}}=0$ with $a_{\sigma_{i}}\neq 0$ for all $i$. For any $\alpha \in \mathbb{E},$ we have $\sum_{i=1}^{k}a_{\sigma_{i}}z_{\sigma_{i}}\alpha=0$, which gives $\sum_{i=1}^{k}\sigma_{i}^{-1}(\alpha)a_{\sigma_{i}}z_{\sigma_{i}}=0$. Consequently, $\sum_{i=1}^{k}\sigma_{i}^{-1}(\alpha)a_{\sigma_{i}}z_{\sigma_{i}}-\sigma_{1}^{-1}(\alpha)\sum_{i=1}^{k}a_{\sigma_{i}}z_{\sigma_{i}}=0,$ i.e. $\sum_{i=2}^{k}(\sigma_{i}^{-1}(\alpha)-\sigma_{1}^{-1}(\alpha))a_{\sigma_{i}}z_{\sigma_{i}}=0,$ a contradiction to $k$ being the least. As $\operatorname{dim}_{\mathbb{E}}(\C_{\mathcal{A}}(\mathcal{B}))= |\G(\mathbb{E}/\mathbb{F})|,$ it follows that $\{z_{\sigma} ~|~ \sigma \in \G(\mathbb{E}/\mathbb{F})\}$ is a basis for the $\mathbb{E}$-vector space $\C_{\mathcal{A}}(\mathcal{B})$. If $\tau: \G(\mathbb{E}/\mathbb{F}) \times \G(\mathbb{E}/\mathbb{F}) \rightarrow \mathcal{U}(\mathbb{E})$ is defined as \begin{equation}\label{neweq2}
\tau(\sigma,\sigma') \mapsto z_{\sigma\sigma'}^{-1} z_{\sigma} z_{\sigma'},
\end{equation} then the associativity of multiplication in $\operatorname{Cen}_{\mathcal{A}}(\mathcal{B})$ implies that $\tau$ is a factor set and $\C_{\mathcal{A}}(\mathcal{B})$ is isomorphic to the classical crossed product $(\mathbb{E}/\mathbb{F},\tau).$ }
\begin{lemma}{\label{l2}}{ Let $(H,K)$  be a generalized strong Shoda pair of $G$, $\lambda$ a linear character of $H$ with kernel $K$. Let $w$ be a normal element of the extension $\mathbb{E}/\mathbb{F}$. Let $\mathcal{G} = \{ \sigma_{1},\sigma_{2},\ldots,\sigma_{\mathtt{k}} \}$ be the Galois group of $\mathbb{E}$ over $\mathbb{F}$ with $\sigma_{1}$ the identity automorphism of $\mathbb{E}.$ If the Schur index of $\mathbb{Q}Ge_{\Q}(\lambda^G)$ is 1, then there exist units $\{\mathbf{z}_{1}, \mathbf{z}_{2}, \ldots , \mathbf{z}_{\mathtt{k}}\}$ in $\operatorname{Cen}_{\mathcal{A}}(\mathcal{B})$ which forms a basis of $\operatorname{Cen}_{\mathcal{A}}(\mathcal{B})$ as a vector space over $\mathbb{E}$ such that the following hold:}
\begin{enumerate}[(i)]
\item  {there exist $\alpha_{1},~\alpha_{2},\ldots, ~\alpha_{\mathtt{k}}$ in $\mathbb{E}$ satisfying the following:\begin{equation}\label{neweq3}
	\begin{array}{lll}\noindent
	\sum_{i=1}^{\mathtt{k}} \alpha_i~\sigma_{i}(w)  &=& \sum_{i=1}^{\mathtt{k}}  \sigma_{i}(w)\\
	\sum_{i=1}^{\mathtt{k}} \alpha_i~\sigma_{i}(\sigma_{j}(w) )&=& w- \sigma_{j}(w),~ 2 \leq j \leq \mathtt{k}
	\end{array}
	\end{equation} and with this choice of $\alpha_i'$s, $\alpha = \sum_{i=1}^{\mathtt{k}} \alpha_{i}\mathbf{z}_{i}$ is invertible. 
	\item  $\{ \mathbf{z}_{i}^{-1} \alpha^{-1} \widehat{E}\, \alpha\, \mathbf{z}_{i}~|~ 1 \leq i \leq {\mathtt{k}}\}$ is a complete set of orthogonal primitive idempotents of $\operatorname{Cen}_{\mathcal{A}}(\mathcal{B})$, where $\widehat{E}=\frac{1}{\mathtt{k}} \sum_{i=1}^{\mathtt{k}} \mathbf{z}_{i}$ and $\alpha$ is as defined in (i).}
\end{enumerate}  	
\end{lemma}
\noindent {\bf Proof.} {For the existence of $\alpha_{i}'$s, it is enough to prove that the ${\mathtt{k}}\times {\mathtt{k}}$ matrix whose $(i,j)$-th entry is $\sigma_{i}(\sigma_{j}(w))$ is invertible. If it is not invertible, then there exist $v_{1},v_{2},\ldots,v_{\mathtt{k}} \in \mathbb{E}$ not all zero such that \begin{equation}\label{neweq12}
\sum_{i=1}^{\mathtt{k}}v_{i}~\sigma_{j}(\sigma_{i}(w))=0,~1\leq j \leq {\mathtt{k}}.
\end{equation}For simplicity, we may assume that $v_{1}\neq 0.$ By multiplying with $v_{1}^{-1}$ if necessary, we may assume that $v_{1}=1.$ For $1\leq j\leq {\mathtt{k}},$ apply $\sigma_{j}^{-1}$ to the equation (\ref{neweq12}) and add over all $j$, we get $\sum_{i=1}^{\mathtt{k}}tr_{\mathbb{E}/\mathbb{F}}(v_{i})\sigma_{i}(w)=0,$ where $tr_{\mathbb{E}/\mathbb{F}}(v_i)$ is the trace of $v_i$. This is not possible because $\sigma_{i}(w),~1\leq i \leq {\mathtt{k}},$ are linearly independent over $\mathbb{F}$ and $tr_{\mathbb{E}/\mathbb{F}}(v_{1})={\mathtt{k}}\neq 0.$\para \noindent As said just before the lemma, for every $i$, there exist invertible $\mathbf{z}_{i}'$s in $\operatorname{Cen}_{\mathcal{A}}(\mathcal{B})$ such that $\sigma_{i}(\alpha)=\mathbf{z}_{i}^{-1}\alpha\mathbf{z}_{i}$, for all $\alpha \in \mathbb{E}$. Since the Schur index of $\Q Ge_{\Q}(\lambda^G)$ is 1, in view of (\cite{her}, Lemma 4.4.1), the twisting $\tau$ in equation (\ref{neweq2}) is equivalent to the trivial twisting. Hence, replacing each $\mathbf{z}_{i}$ with a suitable $\mathbb{E}$-multiple, we may assume that \begin{equation}{\label{eq2}}
\mathbf{z}_{i}\mathbf{z}_{j}=\mathbf{z}_{l},~\mbox{if}~ \sigma_{i} \circ \sigma_{j}=\sigma_{l},~1 \leq i,j \leq {\mathtt{k}}.
\end{equation}\noindent Let $ B=\{\sigma_{1}(w),\sigma_{2}(w),\ldots,\sigma_{\mathtt{k}}(w)\}$. Consider the map $\psi: \operatorname{Cen}_{\mathcal{A}}(\mathcal{B}) \rightarrow M_{\mathtt{k}}(\mathbb{F})$ given by $$\sum_{i=1}^{\mathtt{k}} a_{i}\mathbf{z}_{i} \mapsto [\sum_{i=1}^{\mathtt{k}} l_{a_{i}} \circ \sigma_{i}]_{B},$$ where $a_{i} \in \mathbb{E},~l_{a_{i}}:\mathbb{E} \rightarrow \mathbb{E}$ is given by $l_{a_{i}}(x)=a_{i}x,~1 \leq i \leq {\mathtt{k}}$ and $[\sum_{i=1}^{\mathtt{k}} l_{a_{i}} \circ \sigma_{i}]_{B}$ is the matrix of the $\mathbb{F}$-endomorphism $\sum_{i=1}^{\mathtt{k}} l_{a_{i}} \circ \sigma_{i}$ of $\mathbb{E}$ relative to the basis $B$. Since the Schur index of $\operatorname{Cen}_{\mathcal{A}}(\mathcal{B})$ is 1, it may be checked that $\psi$ is an $\mathbb{F}$-algebra homomorphism. Also, $\psi$ is one-to-one, because $\psi \neq 0$ and $\operatorname{Cen}_{\mathcal{A}}(\mathcal{B})$ is a simple ring. Since the dimension of $M_{\mathtt{k}}(\mathbb{F})$ and $\operatorname{Cen}_{\mathcal{A}}(\mathcal{B})$ over $\mathbb{F}$ are same, it turns out that $\psi$ is surjective and hence it is an isomorphism.\para \noindent 
From the definition of $\psi$, we can see that
\begin{equation}\label{neweq11}
	\psi(\widehat{E}) = [\frac{1}{\mathtt{k}}\sum_{i=1}^{\mathtt{k}} \sigma_{i}]_{B}= \frac{1}{\mathtt{k}}\begin{pmatrix}
		1 & 1 & \cdots & 1 & 1\\
		1 & 1 & \cdots & 1 & 1\\
		\vdots & \vdots & \ddots & \vdots & \vdots \\
		1 & 1 & \cdots & 1 & 1\\
		1 & 1 & \cdots & 1 & 1\\
	\end{pmatrix}.
\end{equation}
Also, if $\alpha_{i}'$s satisfy the system of equations (\ref{neweq3}),  then we can  easily compute the matrix  $[\sum_{i=1}^{\mathtt{k}} l_{\alpha_{i}} \circ \sigma_{i}]_{B}$ and see that \begin{equation}\label{neweq10}
\psi( \sum_{i=1}^{\mathtt{k}}\alpha_i  \mathbf{z}_{i})= [\sum_{i=1}^{\mathtt{k}} l_{\alpha_{i}} \circ \sigma_{i}]_{B}= \begin{pmatrix}
	1 & 1 & 1 & \cdots & 1 & 1\\
	1 & -1 & 0 & \cdots & 0 & 0\\
	1 & 0 & -1 & \cdots & 0 & 0\\
	\vdots & \vdots & \vdots & \ddots & \vdots & \vdots \\
	1 & 0 & 0 & \cdots & -1 & 0\\
	1 & 0 & 0 & \cdots & 0 & -1
\end{pmatrix},
\end{equation}
 which is an invertible matrix. Hence  $\alpha= \sum_{i=1}^{\mathtt{k}}\alpha_i  \mathbf{z}_{i}$ is invertible, which proves $(i).$ To prove $(ii),$ we first  observe that $\hat{E}=\frac{1}{\mathtt{k}} \sum_{i=1}^{n} \mathbf{z}_{i}$ is an idempotent as $\mathbf{z}_{i}'$s satisfy equation (\ref{eq2}). In view of equations (\ref{neweq11}) and (\ref{neweq10}), it turns out that\begin{equation}{\label{eq3}}
 \psi(\alpha^{-1}\widehat{E}\alpha) = \begin{pmatrix}
 1 & 0 & \cdots & 0 & 0\\
 0 & 0 & \cdots & 0 & 0\\
 \vdots & \vdots & \ddots & \vdots & \vdots \\
 0 & 0 & \cdots & 0 & 0\\
 0 & 0 & \cdots & 0 & 0\\
 \end{pmatrix} {\rm ~is~a~primitive~idempotent.}
 \end{equation}
\noindent Consequently, for any  $1 \leq i \leq {\mathtt{k}}$, $\mathbf{z}_{i}^{-1}\alpha^{-1}\widehat{E}\alpha \mathbf{z}_{i}$ is also a  primitive idempotent of $\operatorname{Cen}_{\mathcal{A}}(\mathcal{B})$. It now remains to prove that $\mathbf{z}_{i}^{-1}\alpha^{-1}\widehat{E}\alpha \mathbf{z}_{i}$, $1 \leq i \leq {\mathtt{k}}$, are all orthogonal. Consider $\beta = \mathbf{z}_{i}^{-1} \alpha^{-1} \widehat{E} \alpha \mathbf{z}_{i}$ and $\gamma = \mathbf{z}_{j}^{-1} \alpha^{-1} \widehat{E} \alpha \mathbf{z}_{j},~1 \leq i,j \leq {\mathtt{k}},~i \neq j.$ We have $\beta\gamma =\mathbf{z}_{i}^{-1} \alpha^{-1} \widehat{E} \alpha \mathbf{z}_{i}\mathbf{z}_{j}^{-1} \alpha^{-1} \widehat{E}\alpha \mathbf{z}_{j}.$ By equation (\ref{eq2}), $\mathbf{z}_{i}\mathbf{z}_{j}^{-1}=\mathbf{z}_{l}^{-1}$ for some $l \neq 1,$ as $i \neq j.$ Therefore, to prove $\beta\gamma =0$, it is enough to prove that $\psi(\alpha^{-1} \widehat{E} \alpha \mathbf{z}_{l}^{-1} \alpha^{-1} \widehat{E} \alpha)=0.$  As $\psi(\alpha^{-1} \widehat{E} \alpha \mathbf{z}_{l}^{-1} \alpha^{-1} \widehat{E} \alpha)=\psi(\alpha^{-1} \widehat{E} \alpha) \psi(\mathbf{z}_{l}^{-1}) \psi(\alpha^{-1} \widehat{E} \alpha)$, in view of equation (\ref{eq3}), it follows that $\psi(\alpha^{-1} \widehat{E} \alpha \mathbf{z}_{l}^{-1} \alpha^{-1} \widehat{E} \alpha)$ is the $(1,1)$ entry of $\psi(\mathbf{z}_{l}^{-1})$. Observe that the $(1,1)$ entry of $\psi(\mathbf{z}_{l}^{-1}) = [\sigma_{l}^{-1}]_{B}$ is zero, as $l \neq 1$. Hence, $\psi(\alpha^{-1} \widehat{E} \alpha \mathbf{z}_{l}^{-1} \alpha^{-1} \widehat{E} \alpha) = 0$, and thus $\beta$ and $\gamma$ are orthogonal. This proves the lemma. \qed
\begin{remark}{\label{remark2}}
In the statement of Lemma~\ref{l2}, the assumption that  the Schur index of $\Q Ge_{\Q}(\lambda^G)$ is 1 was necessary in order to prove that $\widehat{E}$ is an idempotent. The cases with the Schur index at least 2 should be handled using a different approach.
\end{remark}
\noindent Now we are ready to state the main result, which is a generalisation of (\cite{JOdRVG2013}, Theorem 4.1), in which we describe a complete set of orthogonal primitive idempotents in a simple component of a rational group algebra with Schur index 1  given by a generalized strong Shoda pair.
\begin{theorem}{\label{thm3}} Let $(H,K)$  be a generalized strong Shoda pair of $G$, {$\lambda$ a linear character of $H$ with kernel $K$} and $H=H_{0}\leq H_{1} \leq \cdots \leq H_{n}=G$ a strong inductive chain from $H$ to $G$. Let $T$, $\widehat{E}$, $\alpha$ and $\mathbf{z}_i'$s be as given by Lemmas~\ref{l1} and \ref{l2}. If the Schur index of $\Q Ge_{\Q}(\lambda^G)$ is 1, then 
	$$\{t^{-1}  \mathbf{z}_{i}^{-1} \alpha^{-1} \widehat{E}\, \varepsilon(H,K)\,  \alpha\, \mathbf{z}_{i}\, t ~|~t \in T,~1 \leq i \leq {\mathtt{k}}\}$$ is a complete set of orthogonal primitive idempotents of $\Q Ge_{\Q}(\lambda^G)$.  \end{theorem}
\noindent\textbf{Proof.} In (\cite{BG}, Section 3), it is proved that the map $\mathcal{B} \otimes_{\mathbb{F}} \operatorname{Cen}_{\mathcal{A}}(\mathcal{B}) \to \Q Ge_{\Q}(\lambda^G)$ given by $x \otimes y \mapsto xy$, for $x \in \mathcal{B}$ and $y \in \operatorname{Cen}_{\mathcal{A}}(\mathcal{B})$ is an $\mathbb{F}$-algebra isomorphism. Therefore, if $\mathcal{E}$ and $\mathcal{E'}$ denote a complete set of orthogonal primitive idempotents of $\mathcal{B}$ and $\operatorname{Cen}_{\mathcal{A}}(\mathcal{B})$ respectively, then $\{ xy ~|~ x \in \mathcal{E},~y \in \mathcal{E'}\}$ is a complete set of orthogonal primitive idempotents of $\Q Ge_{\Q}(\lambda^G)$. In view of Lemmas \ref{l1} and \ref{l2}, $\mathcal{E}=\{t^{-1}\varepsilon(H,K)t~|~t \in T\}$ and $\mathcal{E'}=\{\mathbf{z}_{i}^{-1}\alpha^{-1}\widehat{E}\alpha \mathbf{z}_{i}~|~1 \leq i \leq {\mathtt{k}}\}$ are complete sets of primitive idempotents of $\mathcal{B}$ and $\operatorname{Cen}_{\mathcal{A}}(\mathcal{B})$ respectively. Consider $x = t^{-1}\varepsilon(H,K)t \in \mathcal{E},$ where $t \in T$ and $y = \mathbf{z}_{i}^{-1} \alpha^{-1} \widehat{E}\alpha \mathbf{z}_{i}\in \mathcal{E'}$, for some $1 \leq i \leq {\mathtt{k}}$. Then $$xy=t^{-1}\varepsilon(H,K)t\mathbf{z}_{i}^{-1} \alpha^{-1} \widehat{E}\alpha \mathbf{z}_{i}.$${Using the explicit knowledge of $\Theta_{n-1}$ in terms of the isomorphism $\theta_{j},~0\leq j\leq n-1$ given in equations (\ref{eq1}), (\ref{eq -1}) and (\ref{eq 0}), we see that for any $t \in T,~\Theta_{n-1}(\varepsilon(H,K)t)$ equals a $k \times k$ matrix with precisely one of its entry in the first column $\varepsilon(H,K)$ and all other entries zero. Since $\varepsilon(H,K) \in F,$ we obtain that $\varepsilon(H,K)t\in \mathcal{B}$. This gives in particular for $t=1$ that $\varepsilon(H,K) \in \mathcal{B}.$ As $y \in \operatorname{Cen}_{\mathcal{A}}(\mathcal{B})$, we get that $y$ commutes with $\varepsilon(H,K)t$, and hence $xy= t^{-1}  \mathbf{z}_{i}^{-1} \alpha^{-1} \widehat{E}\alpha \mathbf{z}_{i}\varepsilon(H,K) t.$ Further, $\mathbf{z}_{i}$ and $\alpha$ commute with $\varepsilon(H,K)$ because $\varepsilon(H,K)\in \mathcal{B}$ and $z_{i},~\alpha \in \operatorname{Cen}_{\mathcal{A}}(\mathcal{B})$. So, we have $xy= t^{-1}  \mathbf{z}_{i}^{-1} \alpha^{-1} \widehat{E}\varepsilon(H,K)\alpha \mathbf{z}_{i} t.$  Hence, the theorem follows. }\qed
\para \noindent 
Working analogous to Theorem~\ref{thm3}, we also obtain the matrix units of $\Q Ge_{\Q}(\lambda^G)$ as described in the following theorem which is a generalisation of  (\cite{JOdRVG2013}, Corollary 4.2).
\begin{theorem}{\label{thm4}}  Let $(H,K)$  be  a generalized strong Shoda pair  of $G$, {$\lambda$ a linear character of $H$ with kernel $K$} and $H=H_{0}\leq H_{1} \leq \cdots \leq H_{n}=G$ a strong inductive chain from $H$ to $G$. If the Schur index of $\Q Ge_{\Q}(\lambda^G)$ is 1, then  $$\{ t^{-1}  \mathbf{z}_{i}^{-1} \alpha^{-1} \widehat{E}\,\varepsilon(H,K)\,  \alpha\, \mathbf{z}_{i'}\, t'~|~t, t' \in T, 1 \leq i, i' \leq {\mathtt{k}}\},$$ is a complete set of matrix units of $\Q Ge_{\Q}(\lambda^G)$, where $T$, $\widehat{E}$, $\alpha$  and $\mathbf{z}_i'$s are as given by Lemmas~\ref{l1} and~\ref{l2}. 
\end{theorem}
\section{Units}
In this section, we will show that the knowledge of a complete set of primitive idempotents in all the  simple components of the rational group algebra $\mathbb{Q}G$, for $G$ a generalized strongly monomial group, allows us to describe  a large subgroup of  $\mathcal{U}(\mathbb{Z}G)$, i.e., a subgroup  of  finite index. We will provide two subgroups of  $\mathcal{U}(\mathbb{Z}G)$ which correspond with upper triangular matrices and lower triangular matrices in the simple components. These two subgroups together with any large subgroup of  $\mathcal{Z}(\mathcal{U}(\mathbb{Z}G))$ will generate a subgroup of finite index in  $\mathcal{U}(\mathbb{Z}G)$. 
\para \noindent For $G$ a finite group and $e$ a primitive central idempotent of the rational group algebra $\mathbb{Q} G$, the simple algebra $\mathbb{Q} Ge$ can be identified with the matrix algebra $M_n(D)$ over a division algebra $D$. Let $\mathcal{O}$ be an order in $D$. For any subset $Q$ of $\mathcal{O}$, denote by $E(Q)$ the subgroup of $SL_n(\mathcal{O})$ generated by all $Q$-elementary matrices, i.e.,  $$E(Q)=\left\langle I_n+qE_{ij} \mid q\in Q, 1\leq i,j \leq n, i\neq j, E_{ij} \text{ a matrix unit} \right\rangle.$$ The following is known (see \cite{JL1}, Theorem 2.2).
\begin{theorem}{\rm (Bass-Vaser\u ste\u\i n-Liehl-Venkataramana)} \label{thm5}
	If $Q$ is an ideal in $\mathcal{O}$, the index of  $E(Q)$ in $SL_n(\mathcal{O})$ is finite, provided $n\geq 3$ , or  $n=2$ and  $\mathcal{U}(\mathcal{O})$ is infinite. 
\end{theorem}
\noindent We will now recall a result of the first author with Kaur (\cite{BK19}, Theorem 4), where they have provided a large subgroup  of  $\mathcal{Z}(\mathcal{U}(\mathbb{Z}G))$,  when $G$ is a generalized strongly monomial group. This large subgroup is  contained in the group generated by generalized Bass units of $\mathbb{Z}G$. Below is a quick recall of generalized Bass units in $\mathbb{Z}G$. 
\para\noindent Given $g \in G$ and $k,m$ positive integers such that $k^{m}\equiv 1\operatorname{mod}|g|$, where $|g|$ is the order of $g$, the following is a unit of $\mathbb{Z}G$: $$u_{k,m}(g)=(1+g+\cdots +g^{k-1})^{m}+\frac{1-k^{m}}{|g|}(1+g+\cdots +g^{|g|-1}).$$ The units of this form are called {\it Bass cyclic units} based on $g$ with parameters $k$ and $m$ and were introduced by Bass \cite{Bas66}. When $M$ is a normal subgroup of $G$, then $$u_{k,m}(1-\widehat{M}+g\widehat{M})=1-\widehat{M}+u_{k,m}(g)\widehat{M}$$ is an invertible element of $\mathbb{Z}G(1-\widehat{M})+\mathbb{Z}G\widehat{M}$. As this is an order in $\mathbb{Q}G$, for each element $b=u_{k,m}(1-\widehat{M}+g\widehat{M}),$ there is a positive integer $n$ such that $b^{n} \in \mathcal{U}(\mathbb{Z}G)$. Let $n_{b}$ denote the minimal positive integer satisfying this condition. The element $$u_{k,m}(1-\widehat{M}+g\widehat{M})^{n_{b}}=1-\widehat{M}+u_{k,mn_{b}}(g)\widehat{M}$$ is called the {\it generalized Bass unit} based on $g$ and $M$ with parameters $k$ and $m$. These units were introduced in \cite{JPS96}. Note that if $M$ is such that $G' \subset M$, then the generalized Bass units are central units in $\mathbb{Z}G$.
\begin{theorem}{\rm (\cite{BK22}, Theorem 4)}\label{thm6}
	Let $G$ be a generalized strongly monomial group. Then the group generated by generalized Bass units of $\mathbb{Z}G$ contains a large subgroup of $\mathcal{Z}(\mathcal{U}(\mathbb{Z}G))$.
\end{theorem}
\noindent We can now use matrix units for $G$ obtained in Theorem \ref{thm4} to provide units which generate a subgroup of finite index in $\mathcal{U}(\mathbb{Z}G)$. It may be recalled from (\cite{JdR2016}, Definition 11.2.2) that a simple finite dimensional rational algebra is said to be exceptional if it is one of the following types:\begin{enumerate}[(i)]
\item a non-commutative division algebra different from a totally definite quaternion algebra,
\item $M_{2}(\mathbb{Q}),$
\item $M_{2}(F)$ with $F$ a quadratic imaginary extension of $\mathbb{Q},$
\item $M_{2}(\frac{a,b}{\mathbb{Q}})$ with $a$ and $b$ negative integers, {i.e., $(\frac{a,b}{\mathbb{Q}})$ is a totally definite quaternion algebra with center $\mathbb{Q}$.}
\end{enumerate} We say that a semisimple finite dimensional algebra has no exceptional components if all its Wedderburn components are not exceptional. In the following theorem, we will make use of Notation \ref{not2}.
{\begin{theorem}\label{thm7} Let $(H,K)$ be a generalized strong Shoda pair of $G$ such that the corresponding simple component of $\Q G$ has Schur index 1 and is not exceptional. Let $\lambda$ be a linear character of $H$ with kernel $K$ and let $e=e_{\Q}(\lambda^G)$. Let $X_e$ be a complete set of matrix units of the simple component $\mathbb{Q}Ge$. Let $w$ be an algebraic integer and a normal element of the extension $E$ over $\mathbb{Q}$. Let $\mathcal{T}$ be a transversal of $\mathcal{G} = \operatorname{Gal}(E/F)$ in $\operatorname{Gal}(E/\mathbb{Q})$. Let $\mathtt{B}=\{\sum_{\sigma\in \mathcal{G}}\sigma(\tau(w))~|~\tau \in \mathcal{T}\}$ and let $c$ be an integer  such  that $c\,\beta\,\mathtt{e}_{ij}\in \mathbb{Z}G$ for all $\mathtt{e}_{ij}\in X_e$ and $\beta \in \mathtt{B}$. Let 
	$$V_e^+=\left\langle  1+ c\,\beta\,\mathtt{e}_{ij} \mid  \mathtt{e}_{ij}\ \in X_e, i> j, \, \beta \in \mathtt{B}\right\rangle,$$
	$$V_e^-=\left\langle  1+ c\, \beta\, \mathtt{e}_{ij} \mid  \mathtt{e}_{ij}\ \in X_e, i<j, \beta \in \mathtt{B}\right\rangle.$$
 Then \begin{enumerate}[(i)]
 \item $\langle V_{e}^{+},V_{e}^{-}\rangle$ is a finitely generated subgroup of $\mathcal{U}(\mathbb{Z}G)$;
 \item $ V_{e}^{+}$ and $V_{e}^{-}$ are nilpotent, and $\langle V_{e}^{+},V_{e}^{-}\rangle$ maps to a subgroup of finite index in $1-e+SL_{[G:H]}(\mathcal{O}_{e})$ in the Wedderburn decomposition of $\Q G$, where $\mathcal{O}_{e}=(\mathbb{Z}H\varepsilon)^{\mathcal{G}}$ is the fixed subring of  $\mathbb{Z}H\varepsilon$ under the action of $\mathcal{G}$.
 \end{enumerate}
\end{theorem}}\noindent\textbf{Proof.} It is proved in the proof of Theorem 5 of \cite{BK22} that the center of $\mathbb{Q}Ge$ is isomorphic to $E^\mathcal{G}$. As $w$ is a normal element of the extension $E$ over $\mathbb{Q}$, it turns out that $\mathtt{B}$ is a $\mathbb{Q}$-basis of  $E^\mathcal{G}$ and hence $F=E^\mathcal{G}=\mathbb{Q}(\mathtt{B})$, the field obtained by adjoining the elements of $\mathtt{B}$ to $\mathbb{Q}$. We have that $\Q Ge\cong M_{k}(E/F,\tau).$ Since the Schur index of  $(E/F,\tau)$ is $1$ and the dimension of $(E/F,\tau)$ over $F$ is $\mathtt{k}^2$, it turns out that $(E/F,\tau)\cong M_{\mathtt{k}}(F)$. Consequently, $\Q Ge\cong$ $M_{k\mathtt{k}}(F)$ =$M_{[G:H]}(F)=$ $M_{[G:H]}(\mathbb{Q}(\mathtt{B})).$ \vspace{.2cm}\\ \noindent  $(i)$ To prove that {$V_e^+$ and $V_e^-$ are subgroups of} $\mathcal{U}(\mathbb{Z}G)$, it is enough to show that the generators of $V_e^+ $ and $V_e^-$ belong to $\mathcal{U}(\mathbb{Z}G)$. Consider a generator of $V_e^+ $ or $V_e^-$, which is of the form $1+c\, \beta\,  \mathtt{e}_{ij}$ for some $\beta =\sum_{\sigma  \in \mathcal {G}}\sigma(\tau(w)) \in \mathtt{B}$ and  $i \neq  j$. By the choice of $c$, we already have that $1+c\, \beta\,  \mathtt{e}_{ij} \in \mathbb{Z}G.$ Observe that the product of $1+c\, \beta\,  \mathtt{e}_{ij}$  and $1-c\, \beta\,  \mathtt{e}_{ij}$ is 1, as $i \neq j$. Consequently, $1+c\, \beta\,  \mathtt{e}_{ij}\in \mathcal{U}(\mathbb{Z}G)$. Hence $(i)$ is proved. \vspace{.2cm}\\ \noindent  $(ii)$ Consider the decomposition $\mathbb{Q}G = \mathbb{Q}Ge \bigoplus \mathbb{Q}G(1-e)$. For any $\beta \in  \mathtt{B}$ and $i \neq j$, $1+c\,\beta\,\mathtt{e}_{ij}$ goes to $1-e$ when projected to the summand $\mathbb{Q}G (1-e)$, while in the summand $\mathbb{Q}Ge$ it is projected to the elementary matrix which has all $1$ on the diagonal and other entries are zero except at the $(i,j)$ place, where it is  $c\,\beta$. Therefore, in the above decomposition $\langle V_{e}^{+},V_{e}^{-}\rangle$ maps to $1-e+E(c\mathbb{Z}[\mathtt{B}])$, where $\mathbb{Z}[\mathtt{B}]$ is a $\mathbb{Z}$-module generated by the elements of $\mathtt{B}$. As $w$ is an algebraic integer, the elements of $\mathbb{Z}[\mathtt{B}]$ are integral over $\mathbb{Q}$ and hence $\mathbb{Z}[\mathtt{B}]\subseteq \mathbb{Z}H\varepsilon,$ the ring of integers of $E$. Also, the elements of $\mathtt{B}$ are fixed by $\mathcal{G},$ therefore $$\mathbb{Z}[\mathtt{B}]\subseteq (\mathbb{Z}H\varepsilon)^{\mathcal{G}}.$$ Next choose an integer $d$ such that $d(\mathbb{Z}H\varepsilon)^{\mathcal{G}}\subseteq \mathbb{Z}[\mathtt{B}].$ It may be noted that it is possible to choose such a $d$ because $\mathtt{B}$ is a $\mathbb{Q}$-basis of $E^{\mathcal{G}}=\mathbb{Q}[\mathtt{B}]$ and $(\mathbb{Z}H\varepsilon)^{\mathcal{G}}$ is a finitely generated $\mathbb{Z}$-module. This gives   $$cd(\mathbb{Z}H\varepsilon)^{\mathcal{G}}\subseteq c\mathbb{Z}[\mathtt{B}]\subseteq (\mathbb{Z}H\varepsilon)^{\mathcal{G}}.$$  Set $Q=c\,d\,(\mathbb{Z}H\varepsilon)^{\mathcal{G}},$  and note that  $\mathcal{O}_{e}$  is an order in $E^\mathcal{G}= \mathbb{Q}(\mathtt{B}) $ and $Q$ an ideal in $\mathcal{O}_{e}$. Also, in view of the above arguments, we have $$1-e+E(Q) \subseteq 1-e+E(c\mathbb{Z}[\mathtt{B}])\subseteq 1-e+SL_{[G:H]}(\mathcal{O}_{e}).$$ Since by Theorem~\ref{thm5}, the index of $1-e+E(Q)$ in $1-e+SL_{[G:H]}(\mathcal{O}_{e})$ is finite, we obtain that the index of $\left\langle V_e^+, V_e^-\right\rangle$ ($\cong$ $1-e+E(c\mathbb{Z}[\mathtt{B}])$) in $1-e + SL_{[G:H]}(\mathcal{O}_e)$ is finite. Since $V_{e}^{+}$ and $V_{e}^{-}$ correspond respectively to lower and upper triangular matrices in $SL_{[G:H]}(\mathcal{O}_{e})$, the groups $V_{e}^{+}$ and $V_{e}^{-}$ are nilpotent. This proves $(ii)$.  ~\qed \para\noindent
In view of the above theorem and (\cite{JdR2016}, Remark 4.6.10, Proposition 5.5.1), we immediately have the following generalisation of (\cite{JOdRVG2013}, Theorem 5.4).
\begin{theorem}\label{cor} Let $G$ be a generalized strongly monomial group such that the Schur index of each simple component of $\mathbb{Q}G$ is $1$ and none of them is exceptional. Let $\mathtt{E}$ be the set of all the distinct primitive central idempotents of $\Q G$. For $e\in \mathtt{E},$ let $V_{e}^{+}$, $V_{e}^{-}$ be the subgroups of $\mathcal{U}(\mathbb{Z}G)$ as defined in Theorem \ref{thm7} and set $V^+= \prod_{e\in \mathtt{E}} V_e^+,~ V^- = \prod_{e\in \mathtt{E}} V_e^-.$ Let $U$ be the subgroup of $\mathcal{Z}(\mathcal{U}(\mathbb{Z}G))$ of finite index given by Theorem~\ref{thm6}. Then  \begin{enumerate}[(i)]
\item $V^+$ and $V^-$ are finitely generated nilpotent subgroups of $\mathcal{U}(\mathbb{Z}G)$; 
\item the subgroup $\left\langle U, V^+, V^-\right\rangle$ has finite index in  $\mathcal{U}(\mathbb{Z}G)$.			
\end{enumerate}
\end{theorem}	
\section{Example 1:  $G = P \rtimes D_{2^n}$}
\subsection{The presentation of the group}
\para\noindent
We consider {a group} $G$ that is a semidirect product of an extraspecial $p$-group with a dihedral group.
Let $p$ be an odd prime, $p \equiv 1 \operatorname{mod} 4$ and let $P$ be the extraspecial $p$-group of exponent $p$. By (\cite{Sha}, section IV, Theorem 2), the order of $P$ is $p^{2r+1}$ for an integer $r \geq 1$, and it has generators $x, y_{1},\ldots, y_{r},z_{1},\ldots, z_{r}$, which satisfy the following relations: \begin{equation*}
\begin{split}
& x^p=y_{i}^p=z_{i}^p=1,~1 \leq i \leq r;\\&
[x,y_{i}]=[x,z_{i}]=[y_{i},y_{j}]=[z_{i},z_{j}]=1,[y_{i},z_{i}]=x,~1 \leq i,j \leq r;\\&
[y_{i},z_{j}]=1,~1 \leq i,j \leq r,~i \neq j.
\end{split}
\end{equation*}\noindent
Write $p-1=2^{n-1}t,$ where $(2,t)=1.$ Let $k$ be an integer of order $2^{n-1} \operatorname{modulo} p$ and let $q$ be the multiplicative inverse of $k \operatorname{modulo} p.$ Consider the automorphisms $\theta_{1}$ and $\theta_{2}$ of $P$ given by $\theta_{1}(x)=x,~\theta_{1}(y_{i})=y_{i}^{k},~\theta_{1}(z_{i})=z_{i}^q$, $\theta_{2}(x)=x^{-1},~\theta_{2}(y_{i})=z_{i}^{-1},~\theta_{2}(z_{i})=y_{i}^{-1}$ for all $1 \leq i \leq r.$ It can be checked that the order of $\theta_{1}$ is $2^{n-1}$ and that of $\theta_{2}$ is $2$. Also, $\theta_{2}^{-1}\circ \theta_{1} \circ \theta_{2}= \theta_{1}^{-1}$. Thus $\langle \theta_{1},\theta_{2} \rangle$ is a subgroup of $\operatorname{Aut}(P)$ isomorphic to the Dihedral group $D_{2^n}$ of order $2^n.$ 

Let $G = P \rtimes D_{2^n}$ be the semidirect product of $P$ with respect to $\langle \theta_{1},\theta_{2} \rangle$. {Note that $|G|=p^{2r+1}2^n.$} The group $G$ is generated by $x, y_{1},\ldots, y_{r},z_{1},\ldots, z_{r},a,b$ with the following relations:
\begin{equation*}
\begin{split}
&~x^p=y_{i}^p=z_{i}^p=a^{2^{n-1}}=b^2=1,~1 \leq i \leq r;\\& [x,y_{i}]=[x,z_{i}]=[y_{i},y_{j}]=[z_{i},z_{j}]=1,[y_{i},z_{i}]=x,~1 \leq i,j \leq r;\\& [y_{i},z_{j}]=1~ \mbox{if}~i \neq j;\\ & [a,x]=1,[a,y_{i}]=y_{i}^{1-k},[a,z_{i}]=z_{i}^{1-q}, 1 \leq i \leq r;\\ & [b,x]=x^2,[b,y_{i}]=z_{i}y_{i},[b,z_{i}]=y_{i}z_{i},[b,a]=a^{2}, 1 \leq i \leq r.
\end{split}
\end{equation*}

\begin{remark}For $p=5$ and $r=1$, we obtain the group $G$ of order 1000 that can be identified in the GAP library as \verb+SmallGroup(1000,92)+. This group is not strongly monomial using a GAP check with the GAP package Wedderga (see \cite{Wedderga}). We will show in the next section that $G$ is generalized strongly monomial and will compute a complete set of primitive orthogonal idempotents and matrix units in all the simple components of $\Q G$ using Theorems \ref{thm3} and \ref{thm4}.
\end{remark}

\subsection{Wedderburn decomposition}
In the following theorem, we will show that the hypothesis for Theorems \ref{thm3} and \ref{thm4} is satisfied for the above group $G$. More precisely, we will provide the complete Wedderburn decomposition of $\Q G$ and show that $G$ is a generalized strongly monomial group and the Schur index of all the simple components of $\Q G$ is $1$.		
\begin{theorem}{\label{thm8}} The group $G= P \rtimes D_{2^n}$ as defined above is generalized strongly monomial and the Wedderburn decomposition of $\Q G$ is given by:
	\begin{flalign*}
	\Q G \cong &~ \Q^{(4)} \bigoplus_{i=2}^{n-1} M_{2}(\Q(\zeta_{2^i}+\zeta_{2^i}^{-1})) \bigoplus M_{2^{n-1}}(\Q(\zeta_{p}+\zeta_{p}^{-1}))^{(4m_{1})}\\ & \bigoplus  M_{2^n}(\Q(\zeta_{p}+\zeta_{p}^k+\cdots + \zeta_{p}^{k^{2^{n-1}-1}}))^{(m_{1})}\bigoplus M_{2^{n}}(\Q(\zeta_{p}+\zeta_{p}^{-1}))^{(m_{2})}\\ & \bigoplus_{i=0}^{n-1} M_{2p^r}(\Q(\zeta_{2^ip}+\zeta_{2^ip}^{-1})) ,
	\end{flalign*} where  $m_{1}=\frac{p^r-1}{p-1}$, $m_{2}=\frac{(p^{2r}-1)-(2^{n-1}+2)(p^r-1)}{(p-1)2^{n-1}}$, $\zeta_{l}$ is a primitive $l$-th root of unity and $M_{n}(D)^{(m)}$ is a direct sum of $m$ copies of $M_{n}(D)$.
\end{theorem} \noindent
 \noindent {The strategy to prove this theorem is to provide a complete and irredundant set of generalized strong Shoda pairs of $G$ and write the corresponding simple component using (\cite{JdR2016}, Lemma 3.3.2, Theorem 3.5.5) and Theorem \ref{thm2}. \para \noindent
We begin by finding the Shoda pairs $(H,K)$ for the specific subgroups $H=G, ~P, ~\langle P,a\rangle,$ $\langle P,b\rangle, ~\langle P,ba\rangle,$ and ~$\langle x,y_{1},\ldots,y_{r},a \rangle$. \para \noindent Observe that $G'$, the commutator subgroup of $G$, is  $\langle P,a^2 \rangle,$ and $G/G'$ is the Klein's $4$-group generated by $G'a$, $G'b$. Hence the Shoda pairs $(H,K)$ of $G$ with $H=G$ are $(G,G)$, $(G,\langle P,a^2,ab \rangle)$, $(G,\langle P,a \rangle)$, $(G,\langle P,a^2,b \rangle)$. These Shoda pairs contribute to four copies of $\Q$ in the Wedderburn decomposition of $\Q G$. 
\begin{lemma}\label{l4}
If $H=\langle P,a \rangle$, then $(H,K)$ is a Shoda pair of $G$ if and only if $K = \langle P,a^{2^i} \rangle,~2 \leq i \leq n-1.$ Moreover, all of these Shoda pairs are strong Shoda pairs and are inequivalent. Furthermore, for $K = \langle P,a^{2^i} \rangle,~2 \leq i \leq n-1$, $$\Q Ge(G,H,K) \cong M_{2}(\Q(\zeta_{2^i}+\zeta_{2^i}^{-1})).$$
\end{lemma}\noindent \textbf{Proof.} Observe that $H' =P,$ as $[y_{i},z_{i}]=x,~[a,y_{i}]=y_{i}^{-k+1},~[a,z_{i}]=z_{i}^{-q+1}$ for all $1 \leq i \leq r$. If $H/K$ is cyclic, then possibilities of $K$ are $\langle P,a^{2^i} \rangle,~0 \leq i \leq n-1.$ Observe that normalizer in $G$ of all of these choices of $K$ is $G.$ Also note that $H/K$ is a maximal abelian subgroup of $N_{G}(K)/K $ if and only if $i \neq 0,1.$ Therefore, $(H,K)$ is a Shoda pair if and only if $i \neq 0,1$. Furthermore, these are strong Shoda pairs in view of Remark (2.2, (ii)). Also the Shoda pairs $(H,K),$ $K = \langle P,a^{2^i} \rangle,~2 \leq i \leq n-1$ are all inequivalent because no two of these $K's$ are conjugate in $G$ (see Remark (\ref{remark1}, (ii)). For the simple component $\Q Ge(G,H,K)$, we apply (\cite{JdR2016}, Theorem 3.5.5) and see immediately that if $K = \langle P,a^{2^i} \rangle,~2 \leq i \leq n-1,$ then $\Q Ge(G,H,K) \cong M_{2}(\Q(\zeta_{2^i}+\zeta_{2^i}^{-1})).$ \qed
\para \noindent For  $H=\langle P,b \rangle$ and $H=\langle P,ba \rangle$, we prove the following:}
\begin{lemma}{\label{l5}}
	For each of the subgroups $H=\langle P,b \rangle$ and $H=\langle P,ba \rangle$, there are $2\frac{(p^r-1)}{(p-1)}$ subgroups $K$ such that $(H,K)$ is a Shoda pair of $G$. All of these Shoda pairs are strong Shoda pairs, no two of these are equivalent, and each of these realize the simple component $\Q Ge(G,H,K)$ isomorphic to $M_{2^{n-1}}(\Q(\zeta_{p}+\zeta_{p}^{-1}))$ in the Wedderburn decomposition of $\Q G$.
\end{lemma}
\noindent \textbf{Proof.}  We first consider the case when $H= \langle P,b \rangle.$  As  $[b,z_{i}]=y_{i}z_{i}$ and $[y_{i},z_{i}]=x$ for all $1 \leq i \leq r$, it turns out that $H'=\langle x, y_{1}z_{1},y_{2}z_{2},\ldots,y_{r}z_{r} \rangle.$ Let $$X=\{K~|~ K \unlhd H, K \neq H, H/K~ \mbox{is cyclic}\}.$$ Since $H/H'$ is isomorphic to $\mathbb{Z}_{p}^{(r)} \bigoplus \mathbb{Z}_{2},$ the cardinality of $X$ is equal to $2(\frac{p^r-1}{p-1})+1.$\para \noindent Observe that $X$ contains $P$. \para \noindent \underline{Step 1.} $(H,P)$ is not a Shoda pair of $G$.\par \vspace{.125cm}\noindent
This follows because when $H=\langle P,b \rangle,~[H,a] \cap H=P$ and hence the condition ${\rm (ii)}$ of a Shoda pair (see Section \ref{prel}) is not satisfied.\para \noindent
\underline{Step 2.} $(H,K),~K \in X\backslash \{P\}$, is a strong Shoda pair. \par \vspace{.125cm}\noindent  Consider any $K \in X\backslash \{P\}$. We assert that $N_{G}(K)= \langle P,b,a^{2^{n-2}} \rangle$, from which it can be immediately seen that the conditions $(SS1)$ and $(SS2)$ hold. Trivially, $a^{2^{n-2}}$ normalizes $K$. To show that $a^i$ does not normalize $K$, it is enough to proof that 
\begin{equation}{\label{eq4}}
	H' \nsubseteq K^{a^i},~\forall ~1 \leq i < 2^{n-2},~K \in X\backslash \{P\}.
\end{equation}
Suppose on the contrary, $H' \subseteq K^{a^i}$ for some  $1 \leq i < 2^{n-2}.$ For $1 \leq j \leq r$, as $y_{j}z_{j} \in K$, we have $a^{-i}y_{j}z_{j}a^{i} \in K^{a^i},$ i.e.,  $y_{j}^{k^i}z_{j}^{q^i} \in K^{a^i}$. But $y_{j}z_{j} \in H' \subseteq K^{a^i}$ implies that $y_{j}^{q^i}z_{j}^{q^i} \in K^{a^i}.$ Consequently, we get $y_{j}^{k^i-q^i} \in K^{a^i},$ which gives $y_{j} \in K^{a^i}$ because $p \nmid k^i-q^i$ if $1 \leq i < 2^{n-2}.$ Therefore, {$y_{j}^{q^{i}}=a^iy_{j}a^{-i} \in K$} and hence $y_{j} \in K.$ This gives $K=P$, which is not the case and hence the assertion is proved. Let us now show that $(SS3)$ holds too. Let $P_{j}=\langle x,y_{j},z_{j} \rangle.$ Since $y_{j} \notin K,$ we must have $P_{j} \cap K= \langle x,y_{j}z_{j} \rangle.$ Let $\lambda$ be a linear character of $H$ with kernel $K$ and let $\lambda'$ be the restriction of $\lambda$ on $P_{j}$. Observe that $\operatorname{ker}\lambda'=P_{j} \cap K= \langle x,y_{j}z_{j} \rangle.$ By (\cite{BK19}, Lemma 1), $$e_{\Q}(\lambda)e_{\Q}(\lambda')=e_{\Q}(\lambda)=e_{\Q}(\lambda')e_{\Q}(\lambda).$$ For any $1 \leq i < 2^{n-2}$, conjugating the above equation by $a^i$, we have $$e_{\Q}(\lambda)^{a^i}e_{\Q}(\lambda')^{a^i}=e_{\Q}(\lambda)^{a^i}=e_{\Q}(\lambda')^{a^i}e_{\Q}(\lambda)^{a^i}.$$ This gives $e_{\Q}(\lambda)e_{\Q}(\lambda)^{a^i}=e_{\Q}(\lambda)e_{\Q}(\lambda')e_{\Q}(\lambda')^{a^i}e_{\Q}(\lambda)^{a^i}.$ But $\operatorname{ker}\lambda'=\langle x,y_{j}z_{j} \rangle$ and $\operatorname{ker}\lambda'^{a^i}=\langle x,y_{j}^{k^i}z_{j}^{q^i} \rangle$ being distinct, we have $e_{\Q}(\lambda')$ and $e_{\Q}(\lambda')^{a^i}$ are distinct primitive central idempotents of $P_{j}.$ Therefore, $e_{\Q}(\lambda')e_{\Q}(\lambda')^{a^i}=0.$ Consequently,\linebreak $e_{\Q}(\lambda)e_{\Q}(\lambda)^{a^i}=0.$ \para \noindent \underline{Step 3.} The Shoda pairs $(H,K),~K \in X \backslash \{P\}$, are all inequivalent. \par \vspace{.125cm}\noindent Let $K_{1},~K_{2} \in X \backslash \{P\}.$ First assume that both $K_{1}$ and $K_{2}$ are contained in $P$. {Consider any $1\leq i\leq 2^{n-1}$. As $P$ is normal in $G$ and $P\subseteq H,$ we obtain that $P \subseteq H^{a^i},$ which gives $H^{a^i} \cap K_2=K_2.$ Similarly, $K_1 \subseteq P$ implies $K_{1}^{a^i}\subseteq P^{a^i}=P$, which is contained in $H$ and hence $H\cap K_{1}^{a^i}=K_{1}^{a^i}.$} In view of equation (\ref{eq4}), $K_{2} \neq K_{1}^{a^i}$ for all $1 \leq i < 2^{n-2}.$ Also, as $K_{1}^{a^{2^{n-2}}}= K_{1}$ and $K_{1} \neq K_{2}$, it follows that $K_{2} \neq K_{1}^{a^i}$ for all $1 \leq i \leq 2^{n-1}$. Hence, $H^{a^i} \cap K_{2} \neq H \cap K_{1}^{a^i}$ for all $1 \leq i \leq 2^{n-1}$. Therefore, by Remark (\ref{remark1}, (i)), $(H,K_{1})$ and $(H,K_{2})$ are inequivalent.\para \noindent Similarly, we can show that when both $K_{1}$ and $K_{2}$ are not contained in $P$ or one of them is not contained in $P$, then $(H,K_{1})$ and $(H,K_{2})$ are inequivalent. Furthermore, the simple component $\Q Ge(G,H,K),~K \in X \backslash \{P\}$, is isomorphic to $M_{2^{n-1}}(\Q(\zeta_{p}+\zeta_{p}^{-1}))$ follows directly from (\cite{JdR2016}, Theorem 3.5.5). \para \noindent The case when $H=\langle P,ba \rangle$ is similar and exactly the same way one can show that $X'=\{K~|~ K \unlhd H, K \neq H, H/K~ \mbox{cyclic}\}$ has cardinality $2(\frac{p^r-1}{p-1})+1$, $(H,P)$ is not a Shoda pair of $G$, $(H,K),~K \in X' \backslash \{P\}$ are all strong Shoda pair, these are all inequivalent, and the simple component corresponding to any of these is isomorphic to $M_{2^{n-1}}(\Q(\zeta_{p}+\zeta_{p}^{-1}))$. \para \noindent Finally, to complete the proof, it remains to show that if $K \in X$ and $K' \in X'$, then $(\langle P,b \rangle,K)$ and $(\langle P,ba \rangle,K')$ are inequivalent. This follows because for all $1 \leq i \leq 2^{n-1}$,\begin{equation*}
\langle P,b \rangle^{a^i} \cap K'=\left\{
\begin{array}{ll}
K' ,& \mbox{if}~ ba \notin K';\\
K'\backslash \{ba\}, &\mbox{if}~ba \in K';
\end{array}
\right.
\mbox{and}~
\langle P,ba \rangle \cap K^{a^i}=\left\{
\begin{array}{ll}
K^{a^i}, & \mbox{if}~ b \notin K;\\
(K\backslash \{b\})^{a^i}, &\mbox{if}~b \in K.
\end{array}
\right.
\end{equation*}~\qed \para \noindent
In the following lemma, we will find the Shoda pairs $(H,K)$ of $G$ when  $H=P$.
\begin{lemma}{\label{l6}}
	If $H=P$, then there are:
	\begin{enumerate}[(i)]
		\item  $m_{1}$ inequivalent Shoda pairs $(H,K)$ of $G$ with simple component isomorphic to $M_{2^n}(\Q(\zeta_{p}+\zeta_{p}^k+\cdots+\zeta_{p}^{k^{2^{n-1}-1}}))$, where $m_{1}=\frac{p^r-1}{p-1}$. 
		\item $m_{2}$ inequivalent Shoda pairs $(H,K)$ of $G$ with simple component isomorphic to $M_{2^n}(\Q(\zeta_{p}+\zeta_{p}^{-1})),$ where $m_{2}=\frac{(p^{2r}-1)-(2^{n-1}+2)(p^r-1)}{(p-1)2^{n-1}}.$ 
	\end{enumerate}All the Shoda pairs in (i) and (ii) above are strong Shoda pairs.
\end{lemma}
\noindent \textbf{Proof.} Note that $H'=P'=\langle x \rangle.$ Let $
X=\{K~|~K \unlhd H, K \neq H, H/K ~\mbox{is cyclic}\}. $  Since $H/H'$ is isomorphic to $\mathbb{Z}_{p}^{(2r)}$, the cardinality of $X$ is $\frac{p^{2r}-1}{p-1}.$ Let \begin{equation}{\label{eq6}}
Y=\{K~|~K \in X, \langle x,y_{1}^{k^j}z_{1},\ldots,y_{r}^{k^j}z_{r}\rangle \subseteq K {\rm ~for~some~}j,~1 \leq j \leq 2^{n-1}\}.
\end{equation} As $H/\langle x,y_{1}^{k^j}z_{1},\ldots,y_{r}^{k^j}z_{r}\rangle \cong \mathbb{Z}_{p}^{(r)},$ it follows that the cardinality of $Y$ is $2^{n-1}(\frac{p^r-1}{p-1}).$ \par \vspace{.125cm}\noindent Note that $(H,K)$ is not a Shoda pair if $K \in Y$. This is because, for any $j$, $[P,ba^j] \cap P \subseteq K$ if $\langle x,y_{1}^{k^j}z_{1},\cdots,y_{r}^{k^j}z_{r}\rangle\subseteq K$. \para \noindent
We now show that $(H,K),~K \in X \backslash Y$, is a strong Shoda pair of $G$. Let \begin{equation}{\label{eqo}}
Z=\{K~|~ K \in X \backslash Y, \langle x,y_{1}, \ldots, y_{r} \rangle \subseteq K~{\rm or}~ \langle x, z_{1},\ldots , z_{r} \rangle \subseteq K \}.
\end{equation} Observe that the cardinality of $Z$ is $2(\frac{p^r-1}{p-1}).$ \para \noindent
\underline{Case (i).} $K \in Z$. \par \vspace{.125cm}\noindent Observe that in this case $a$ normalizes $K$ and $b$ does not. Therefore, $N_{G}(K)=\langle P,a \rangle$. Also, we have that for some $i$, either $y_{i}$ or $z_{i}$ does not belong to $K$. If $y_{i} \notin K$ (respectively $z_{i} \notin K$), then $a^jK$ can not commute with $y_{i}K$ (respectively $z_{i}K$), for all $1 \leq j \leq 2^{n-1}-1.$ Thus, $H/K$ is a maximal abelian subgroup of $N_{G}(K)/K$. {In view of Remark (\ref{remark3}, (i)), the condition $(SS1)$ and $(SS3)$ hold because $H$ is normal in $G$ and hence $(H,K)$ is a strong Shoda pair of $G$, if $K \in Z$.} Since each $K \in Z$ has precisely two conjugates in $G$, namely $K$ and $K^b$, therefore, we obtain $\frac{p^r-1}{p-1}$ inequivalent Shoda pairs from this list in view of Remark (\ref{remark1}, (ii)). The simple component corresponding to any of these Shoda pairs is isomorphic to $M_{2^{n}}(\Q(\zeta_{p}+\zeta_{p}^k+\cdots +\zeta_{p}^{k^{2^{n-1}-1}}))$ follows from (\cite{JdR2016}, Theorem 3.5.5).\para \noindent
\underline{Case (ii).} $K \notin Z,~K \notin Y.$\par \vspace{.125cm}\noindent {Clearly, $a^{2^{n-2}}$ normalizes $K$, and hence, $\langle P,a^{2^{n-2}}\rangle \subseteq N_{G}(K)$. In this case, there must exist $i_{1},i_{2},\cdots,i_{m},j_{1},\cdots,j_{m'}$ such that $y_{i_{1}}^{\alpha_{1}}y_{i_{2}}^{\alpha_{2}}\cdots y_{i_{m}}^{\alpha_{m}}, z_{i_{1}}^{\beta_{1}}z_{i_{2}}^{\beta_{2}}\cdots z_{i_{m'}}^{\beta_{m'}} \notin K$ but $
y_{i_{1}}^{\alpha_{1}}\cdots y_{i_{m}}^{\alpha_{m}}z_{i_{1}}^{\beta_{1}}\cdots z_{i_{m'}}^{\beta_{m'}} \in K$  for some $\alpha_{i}$ and $\beta_{j}$ non zero. If $a^l$ normalizes $K$ for some $1\leq l<2^{n-2}$, then $ 
a^{-l}y_{i_{1}}^{\alpha_{1}}\cdots y_{i_{m}}^{\alpha_{m}}z_{i_{1}}^{\beta_{1}}\cdots z_{i_{m'}}^{\beta_{m'}}a^l\in K$, i.e.,   $ 
y_{i_{1}}^{\alpha_{1}k^l}\cdots y_{i_{m}}^{\alpha_{m}k^l}z_{i_{1}}^{\beta_{1}q^l}\cdots z_{i_{m'}}^{\beta_{m'}q^l}\in K.$ Also as $(y_{i_{1}}^{\alpha_{1}}\cdots y_{i_{m}}^{\alpha_{m}}z_{i_{1}}^{\beta_{1}}\cdots z_{i_{m'}}^{\beta_{m'}})^{q^l} \in K$,  one gets $y_{i_{1}}^{\alpha_{1}q^l}\cdots y_{i_{m}}^{\alpha_{m}q^l}z_{i_{1}}^{\beta_{1}q^l}\cdots z_{i_{m'}}^{\beta_{m'}q^l}x^{\alpha}\in K$, for some integer $\alpha\geq 1.$ Since $x \in K$,  it yields 
$y_{i_{1}}^{\alpha_{1}q^l}\cdots y_{i_{m}}^{\alpha_{m}q^l}z_{i_{1}}^{\beta_{1}q^l}\cdots z_{i_{m'}}^{\beta_{m'}q^l}\in K$. This gives 
$y_{i_{1}}^{\alpha_{1}(q^l-k^l)}\cdots y_{i_{m}}^{\alpha_{m}(q^l-k^l)}\in K,$ and consequently  $y_{i_{1}}^{\alpha_{1}}\cdots y_{i_{m}}^{\alpha_{m}} \in K$ as $p \nmid q^l-k^l$, which is a contradiction. Hence we have shown that $a^l,~1 \leq l < 2^{n-2}$, can not normalize $K$.\para \noindent Also, we will show that $ba^l$ does not normalize $K$ for any $0\leq l < 2^{n-2}$. We will prove this by contradiction. Suppose $ba^l\in N_{G}(K)$ for some $l,~0 \leq l <2^{n-2}$. We assert that, in this situation, either $y_{i}^{k^l}z_{i}\in K$ for all $i$ or $y_{i}^{-k^l}z_{i}\in K$ for all $i$, which implies $K \in Y,$ not the case. Consider any $1 \leq i \leq n.$ If $i$ is such that $y_{i} \in K$, then $z_{i}^{-1}=(ba^l)^{-1}y_{i}^{k^l}(ba^l)\in K$. Hence both $y_{i}^{k^l}z_{i}$ and $y_{i}^{-k^l}z_{i}$ belong to $K$. If $y_{i}\notin K,$ observe that $H/K=\langle y_{i}K\rangle.$ As $H\unlhd G,$ we get $(ba^l)^{-1}y_{i}^{k^l}(ba^l)K=y_{i}^{m}K$ for some integer $m$. Since $ba^l$ has order 2, $m=\pm k^l.$ Consequently, either $y_{i}^{k^l}z_{i}\in K$ or $y_{i}^{-k^l}z_{i}\in K$. Next assume $i$ and $j$ are such that both $y_{i}$ and $y_{j}$ are not in $K$ but $y_{i}^{k^l}z_{i}\in K$ and $y_{j}^{-k^l}z_{j}\in K$. Arguing as before, we can see that at least one of $(y_{i}y_{j})^{k^l}z_{i}z_{j}$ or $(y_{i}y_{j})^{-k^l}z_{i}z_{j}$ belong to $K$. If $(y_{i}y_{j})^{k^l}z_{i}z_{j}$ (respectively $(y_{i}y_{j})^{-k^l}z_{i}z_{j}$), $y_{i}^{k^l}z_{i}$ and $y_{j}^{-k^l}z_{j}$ belong to $K$, then it turns out that $y_{j}$ (respectively $y_{i}$) belong to $K$, a contradiction. This proves the assertion. \para \noindent Consequently, $N_{G}(K)=\langle P,a^{2^{n-2}} \rangle$. The condition $(SS2)$ holds because if $y_{i} \notin K$ (respectively $z_{i} \notin K$), then $a^{2^{n-2}}K$ can not commute with $y_{i}K$ (respectively $z_{i}K$). Further, $(SS1)$ and $(SS3)$ hold because $H$ is normal in $G$. Also, because $[G:N_{G}(K)]=2^{n-1}$, each $K$ has $2^{n-1}$ conjugates in $G$. Therefore, there are $\frac{(p^{2r}-1)-(2^{n-1}+2)(p^r-1)}{(p-1)2^{n-1}}$  inequivalent Shoda pairs from this list. By using (\cite{JdR2016}, Theorem 3.5.5), the simple component corresponding to any of these Shoda pairs turns out to be isomorphic to $M_{2^{n}}(\Q(\zeta_{p}+\zeta_{p}^{-1})).$ }~\qed \para\noindent
Now we will find the generalized strong Shoda pairs $(H,K)$ of $G$ when  $H=\langle x,y_{1},\ldots,y_{r},a \rangle$. 
\begin{lemma}{\label{l7}}
	 If ${H}=\langle x,y_{1},\ldots,y_{r},a \rangle$, then $(H,K)$ is a Shoda pair of $G$ if and only if $K = \langle y_{1},\ldots,y_{r},a^{2^j} \rangle,~0 \leq j \leq n-1.$ All of these Shoda pairs are generalized strong Shoda pairs and are inequivalent. Furthermore,
	 \begin{enumerate}[(i)]
	 \item  if $H_{1}=\langle P,a\rangle,$ then $H_{0}=H\leq H_{1}\leq H_{2}=G$ is a strong inductive chain from $H$ to $G$.
	 \item if $\lambda$ is a linear character of $H$ with kernel $K$, then $C_{0}=\operatorname{Cen}_{H_{1}}(e_{\Q}(\lambda))=H$ and $C_{1}=\operatorname{Cen}_{G}(e_{\Q}(\lambda^{H_{1}}))=G.$
	 \end{enumerate}  \end{lemma}
\noindent \textbf{Proof.} Observe that $[a,y_{i}]=y_{i}^{-k+1}$ for all $1 \leq i \leq r$ and ${H}/\langle y_{1},\ldots,y_{r} \rangle$ is abelian. Therefore, ${H}'= \langle y_{1}, \ldots ,y_{r} \rangle.$ If $(H,K)$ is a Shoda pair, then ${H}' \subseteq K$ implies the following possibilities of $K$: $$ \langle x,y_{1},\ldots,y_{r},a^{2^j}\rangle \text{ or }\langle y_{1},\ldots,y_{r},a^{2^j}\rangle,~0 \leq j \leq n-1.$$ If $K=\langle x,y_{1},\ldots,y_{r},a^{2^j}\rangle$ for some $0 \leq j \leq  n-1$, then $[H,z_{1}]\cap H =\langle x\rangle$, which is contained in $K$. Hence, the condition (ii) from the definition of a Shoda pair (see Section~\ref{prel}) is violated and $(H,K)$ cannot be a Shoda pair.\para \noindent We will now prove that if $K=\langle y_{1},\ldots,y_{r},a^{2^j} \rangle,~0 \leq j \leq n-1,$ then $(H,K)$ is a generalized strong Shoda pair with strong inductive chain $$H_{0}={H} \leq H_{1}=\langle P,a \rangle \leq G.$$ To prove this, in view of (\cite{BK19}, Theorem 2), it is enough to show that if $\lambda$ is a linear character of $H$ with kernel $K=\langle y_{1},\ldots,y_{r},a^{2^j} \rangle$, then the following holds for $i=0,1$: \begin{enumerate}[(a)]
	\item $H_{i} \unlhd \operatorname{Cen}_{H_{i+1}}(e_{\Q}(\lambda^{H_{i}}));$
	\item the distinct $H_{i+1}$-conjugates of $e_{\Q}(\lambda^{H_{i}})$ are mutually orthogonal;
	\item $(\lambda^{H_{i}})^x=\lambda^{H_{i}}$ for $x \in \operatorname{Cen}_{H_{i+1}}(e_{\Q}(\lambda^{H_{i}}))$ implies that $x \in H_{i}$.
\end{enumerate} For $i=0$, (a) and (c) holds trivially because $\operatorname{Cen}_{H_{1}}(e_{\Q}(\lambda))=N_{H_{1}}(K)={H}$ by (\cite{JdR2016}, Lemma 3.5.1). Let us now prove (b) for $i=0$. Let $\lambda'$ be the restriction of $\lambda$ on $\langle x,y_{1},\ldots ,y_{r} \rangle.$ By ({\cite{BK19}}, Lemma 1), $$e_{\Q}(\lambda)e_{\Q}(\lambda')=e_{\Q}(\lambda)=e_{\Q}(\lambda')e_{\Q}(\lambda).$$ Conjugating the above equation by {$z=z_{1}^{\alpha_{1}}\cdots z_{r}^{\alpha_{r}},$} where $0 \leq \alpha_{i} \leq p-1$, $z \neq 1$, we have $$e_{\Q}(\lambda)^{z}e_{\Q}(\lambda')^{z}=e_{\Q}(\lambda)^{z}=e_{\Q}(\lambda')^{z}e_{\Q}(\lambda)^{z}.$$ This gives \begin{equation}{\label{eq7}}
e_{\Q}(\lambda)e_{\Q}(\lambda)^{z}=e_{\Q}(\lambda)e_{\Q}(\lambda')e_{\Q}(\lambda')^{z}e_{\Q}(\lambda)^{z}.
\end{equation}Now, $\operatorname{ker}\lambda'=\langle y_{1},\ldots ,y_{r} \rangle$ and $\operatorname{ker}\lambda'^z =\langle y_{1}x^{\alpha_{1}},\ldots ,y_{r}x^{\alpha_{r}} \rangle$ being distinct, we have that $e_{\Q}(\lambda')$ and $e_{\Q}(\lambda')^{z}$ are distinct primitive central idempotents of the group algebra $\Q \langle x,y_{1},\ldots ,y_{r} \rangle.$ Therefore, $e_{\Q}(\lambda')e_{\Q}(\lambda')^{z}=0.$ Consequently, equation (\ref{eq7}) implies that $e_{\Q}(\lambda)e_{\Q}(\lambda)^{z}=0.$ This proves (b) for $i=0.$
\para \noindent Suppose $i=1.$ We first show that $\operatorname{Cen}_{G}(e_{\Q}(\lambda^{H_{1}}))=G$. For any $0\leq l_{1},~ \beta_{i},~\gamma_{i}\leq p-1,~1\leq i \leq r$ and $0 \leq l_{2}\leq 2^{n-1}-1$, from (\cite{Iss}, Definition 5.1), we have that
\begin{equation*}\scriptsize
\begin{array}{lll}\noindent
\lambda^{H_{1}}(x^{l_{1}}y_{1}^{\beta_{1}}\cdots y_{r}^{\beta_{r}}z_{1}^{\gamma_{1}}\cdots z_{r}^{\gamma_{r}}a^{l_{2}}) &=&\underset{\stackrel{0 \leq \alpha_{i}\leq p-1}{1\leq i \leq r}}{\sum}\lambda^{\circ}((z_{1}^{\alpha_{1}}\cdots z_{r}^{\alpha_{r}})^{-1}x^{l_{1}}y_{1}^{\beta_{1}}\cdots y_{r}^{\beta_{r}}z_{1}^{\gamma_{1}}\cdots z_{r}^{\gamma_{r}}a^{l_{2}}z_{1}^{\alpha_{1}}\cdots z_{r}^{\alpha_{r}})\\
&=& \underset{\stackrel{0 \leq \alpha_{i}\leq p-1}{1\leq i \leq r}}{\sum}\lambda^{\circ}(x^{l_1+\sum_{i=1}^{r}\alpha_{i}\beta_{i}}y_{1}^{\beta_{1}}\cdots y_{r}^{\beta_{r}}z_{1}^{\gamma_{1}-\alpha_{1}+k^{l_2}\alpha_{1}}\cdots z_{r}^{\gamma_{r}-\alpha_{r}+k^{l_2}\alpha_{r}}a^{l_2}).
\end{array}
\end{equation*} Here, $\lambda^{\circ}(\alpha)$ denotes $\lambda(\alpha)$ if $\alpha \in H$ and $0$ otherwise. On simplifying the above expression, it turns out that\begin{equation}{\label{eq9}}
\lambda^{H_{1}}(x^{l_{1}}y_{1}^{\beta_{1}}\cdots y_{r}^{\beta_{r}}z_{1}^{\gamma_{1}}\cdots z_{r}^{\gamma_{r}}a^{l_{2}}) = \left\{
\begin{array}{ll}
\zeta_p^{l_{3}}\zeta_{2^j}^{l_{2}}, & \mbox{if}~l_{2} \neq 0;\\
p^r \zeta_p^{l_{1}}, &\mbox{if}~ l_{2}=\beta_{i}=\gamma_{i}=0,~1 \leq i \leq r;\\
0,&\mbox{if}~l_{2}=0,~\beta_{i}~\mbox{or}~\gamma_{i}\neq 0~\mbox{for some}~i,
\end{array}
\right.
\end{equation}
where {$l_{3}=l_{1}+\sum_{i=1}^{r}\frac{\beta_{i}\gamma_{i}}{1-k^{l_{2}}}$}. {Similarly, note that} 
\begin{equation}{\label{eq10}}
(\lambda^{H_{1}})^{b}(x^{l_{1}}y_{1}^{\beta_{1}}\cdots y_{r}^{\beta_{r}}z_{1}^{\gamma_{1}}\cdots z_{r}^{\gamma_{r}}a^{l_{2}}) = \left\{
\begin{array}{ll}
\zeta_p^{-l_{3}}\zeta_{2^j}^{-l_{2}}, & \mbox{if}~l_{2} \neq 0;\\
p^r \zeta_p^{-l_{1}}, &\mbox{if}~ l_{2}=\beta_{i}=\gamma_{i}=0,~1 \leq i \leq r;\\
0,&\mbox{if}~l_{2}=0,~\beta_{i}~\mbox{or}~\gamma_{i}\neq 0~\mbox{for some}~i.
\end{array}
\right.
\end{equation}
Therefore, if $\sigma \in \operatorname{Aut}(\Q(\zeta_{2^{n-1}p}))$ such that $\sigma(\zeta_{p})=\zeta_{p}^{-1}$ and $\sigma(\zeta_{2^{n-1}})=\zeta_{2^{n-1}}^{-1}$, then $(\lambda^{H_{1}})^b=\sigma(\lambda^{H_{1}}).$ This implies $(e_{\Q}(\lambda^{H_{1}}))^b=e_{\Q}(\lambda^{H_{1}}),$ i.e., $b \in \operatorname{Cen}_{G}(e_{\Q}(\lambda^{H_{1}})).$ Hence, $\operatorname{Cen}_{G}(e_{\Q}(\lambda^{H_{1}}))=G.$ Now (a) and (b) {hold} trivially because $H_{1} \unlhd G$. We note from equation (\ref{eq9}) and (\ref{eq10}) that $(\lambda^{H_{1}})^{b}(x)=p^r\lambda(x^{-1})$, $\lambda^{H_{1}}(x)=p^r\lambda(x)$. Since $x^2 \notin \operatorname{ker}\lambda,$ $(\lambda^{H_{1}})^b \neq \lambda^{H_{1}}$ and hence the condition (c) holds for $i=1$. We have thus obtained the following {generalized strong} Shoda pairs of $G$:
\begin{equation}{\label{eq8}}
({H},\langle y_{1},\ldots,y_{r},a^{2^j}\rangle),~0 \leq j \leq n-1.
\end{equation}
\para \noindent We will now show that all the generalized strong Shoda pair of $G$ in equation (\ref{eq8}) are inequivalent. Let $0 \leq l,m \leq n-1$ and let $\lambda_{1}$ (respectively $\lambda_{2}$) be a linear character of $H$ with $K=\langle y_{1},\ldots , y_{r},a^{2^l}\rangle$ (respectively $K=\langle y_{1},\ldots , y_{r},a^{2^m} \rangle$). From equations (\ref{eq9}) and (\ref{eq10}), $\lambda_{1}^G(a)=\zeta_{2^l}+\zeta_{2^l}^{-1}$ and $\lambda_{2}^G(a)=\zeta_{2^m}+\zeta_{2^m}^{-1}$. If $l \neq m$, then $\zeta_{2^l}+\zeta_{2^l}^{-1}$ and $\zeta_{2^m}+\zeta_{2^m}^{-1}$ have minimal polynomial over $\Q$ of distinct degree. Therefore, there is no $\sigma \in \operatorname{Aut}(\mathbb{C})$ such that $\lambda^{G}_{1}(a)=\sigma \circ \lambda_{2}^G(a)$. This proves, in view of Remark (\ref{remark1}, (iii)), that all of these generalized strong Shoda pairs are inequivalent. This completes the proof of the lemma. ~\qed 
\begin{lemma}{\label{l8}}
	If $H={\langle x,y_{1},\ldots,y_{r},a \rangle}$ and $K= \langle y_{1},\ldots,y_{r},a^{2^j} \rangle,~0 \leq j \leq n-1,$ then $$\Q Ge_{\Q}(\lambda^G) \cong M_{2p^r}(\Q(\zeta_{2^jp}+\zeta_{2^jp}^{-1})),$$ where $\lambda$ is a linear character of $H$ with kernel $K$. 
\end{lemma}
\noindent \textbf{Proof.} In Lemma~\ref{l7}, we have seen that $(H,K)$ is a generalized strong Shoda pair of $G$ with $H=H_{0}\leq H_{1} \leq G$ as a strong inductive chain, where $H_{1}=\langle P,a \rangle$. Also, we have seen in the same lemma that $C_{0}=\operatorname{Cen}_{H_{1}}(e_{\Q}(\lambda))={H}$ and $C_{1}=\operatorname{Cen}_{G}(e_{\Q}(\lambda^{H_{1}}))=G.$
\para \noindent Let $T=\langle z_{1},z_{2},\cdots,z_{r}\rangle=\{z_{1}^{i_{1}}z_{2}^{i_{2}}\cdots z_{r}^{i_{r}}~|~0 \leq i_{j} \leq p-1,~1 \leq j \leq r \}$. As $C_{0}={H}$ and $H_{1}=\langle P,a\rangle$, $T$ is a transversal of $C_{0}$ in $H_{1}.$ We order the elements of $T$ as follows: first all powers of $z_{1}$ (in ascending order), followed by the power of $z_{1}$ times $z_{2}$,  followed by the power of $z_{1}$ times $z_{2}^2$ and so on. Note that the first element of $T$ is identity. Denote by $\beta_{j}$ the $j$-th element of $T,~1 \leq j \leq p^r.$ We have the following observations:
\begin{enumerate}[(i)]
	\item From (\cite{BK19}, Proposition 2), the map $\theta: \Q H_{1}e_{\Q}(\lambda^{H_{1}}) \rightarrow M_{p^r}(\Q He_{\Q}(\lambda))$ given by $\alpha \mapsto (\alpha_{ij})$ is an isomorphism, where $\alpha_{ij}=e_{\Q}(\lambda)\beta_{j}\alpha\beta_{i}^{-1}e_{\Q}(\lambda).$
	\item From Theorem~\ref{thm1}, $\Q Ge_{\Q}(\lambda^G) \cong \Q H_{1}e_{\Q}(\lambda^{H_{1}})*^{\sigma_{1}}_{\tau_{1}}G/H_{1}$, where the twisting $\tau_{1}$ is trivial and $\sigma_{1}: G/H_{1} \rightarrow \operatorname{Aut}(\Q H_{1}e_{\Q}(\lambda^{H_{1}}))$ sends $bH_{1}$ to $(\sigma_{1})_{b}$. Here $(\sigma_{1})_{b}$ is the automorphism of $\Q H_{1}e_{\Q}(\lambda^{H_{1}})$ given by conjugation by $b$.
	\item From (i) and (ii), $\Q Ge_{\Q}(\lambda^G) \cong M_{p^r}(\Q He_{\Q}(\lambda))*^{\sigma_{1}}_{\tau_{1}}G/H_{1}$.
	\item Theorem~\ref{thm2} implies that $\Q Ge_{\Q}(\lambda^G) \cong M_{p^r}(\Q He_{\Q}(\lambda))*^{\sigma_{1}}_{\tau_{1}}G/H \cong M_{p^r}(\mathbb{E}*_{\tau}\mathcal{G})$, where $\mathbb{E}=\{\operatorname{diag}(\alpha,\ldots,\alpha)_{p^r}~|~\alpha \in \Q He_{\Q}(\lambda)\}$, $\mathbb{F}=\mathcal{Z}(\Q Ge_{\Q}(\lambda^G))$, $\mathcal{G}=\operatorname{Gal}(\mathbb{E}/\mathbb{F})$ is of type $C_{0}/H_{0}$-by-$C_{1}/H_{1}$ and $\tau$ is as defined in equation (\ref{n}). As $C_{0}=H_{0}$, $C_{1}=G$ and $[G:H_{1}]=2,$ it turns out that $\mathcal{G}$ is a cyclic group of order $2$.
	\item Further, it follows from the {proof of Theorem~\ref{thm2} given in (\cite{BG}, section 3)} that there exists $A_{b} \in M_{p^r}(\Q He_{\Q}(\lambda))$, such that $z_{b}=bA_{b}$ belongs to $\operatorname{Cen}_{\mathcal{A}}(\mathcal{B})$, where $\mathcal{B}=M_{p^r}(F)$, $\mathcal{A}=M_{p^r}(\Q He_{\Q}(\lambda))*G/H_{1}$ and $F=\{\alpha \in \Q He_{\Q}(\lambda)~|~ \operatorname{diag}(\alpha,\alpha,\ldots,\alpha)_{p^r} \in \mathbb{F} \}$.
	\item With the choice of $z_{b}$ given above in (v), it turns out that conjugation by $z_{b}$ on $\mathbb{E}$ is an automorphism of $\mathbb{E}/\mathbb{F}$ of order $2$. Denoting this automorphism by $\sigma_{b}$, we have $\mathcal{G}=\langle \sigma_{b} \rangle$ and \begin{equation}{\label{nn}}
	\Q Ge_{\Q}(\lambda^G) \cong M_{p^r}(\mathbb{E}*_{\tau}\mathcal{G})\cong M_{p^r}((\mathbb{E}/\mathbb{F},\sigma_{b},z_{b}^2)).
	\end{equation} \end{enumerate}We now wish to show that $\Q Ge_{\Q}(\lambda^G) \cong M_{2p^r}(\mathbb{F}),$ which can be achieved if we show that 
	\begin{equation}{\label{eq11}}
	z_{b}^2 \in N_{\mathbb{E}/\mathbb{F}}(\alpha)~\mbox{for some}~0 \neq \alpha \in \mathbb{E}.
	\end{equation}  The problem thus boils down to find $A_{b} \in M_{p^r}(\Q H e_{\Q}(\lambda))$ so that $z_{b}=bA_{b} \in \operatorname{Cen}_{\mathcal{A}}(\mathcal{B})$ (which definitely exists by Theorem~\ref{thm2}) and show that $z_{b}^2$ satisfies equation (\ref{eq11}). \para \noindent Suppose $A_{b}=(a_{ij})_{p^r \times p^r},~a_{ij} \in \Q He_{\Q}(\lambda)$ so that $bA_{b} \in \operatorname{Cen}_{\mathcal{A}}(\mathcal{B}).$ Let $B_{k}$ be a $p^r \times p^r$ matrix whose $(1,k)$-th entry is $e_{\Q}(\lambda)$ and all other entries are $0$. Clearly, we have that $B_{k} \in \mathcal{B}$. Hence,
\begin{equation}{\label{eq12}}
A_{b}^{-1}b^{-1}B_{k}bA_{b}=B_{k},~\mbox{for all}~ 1 \leq k \leq p^r.
\end{equation} \noindent For $k=1,$ it gives 
$$A_{b}^{-1}b^{-1}B_{1}bA_{b}=B_{1}.$$ Let us denote $b^{-1}B_{1}b$ by ${X_{1}}$ and let $c_{ij}$ be its $(i,j)$-th entry. From the isomorphism $\theta$, $$c_{ij}=e_{\Q}(\lambda)\beta_{j}b^{-1}e_{\Q}(\lambda)b\beta_{i}^{-1}e_{\Q}(\lambda),$$ which can be seen to be equal to $\frac{e_{\Q}(\lambda)}{p^r}$ using the fact that $y_{i} \in K~\forall ~i$ and group relations {$[b,x]=x^2,[b,y_{i}]=z_{i}y_{i},[b,z_{i}]=y_{i}z_{i},[b,a]=a^{2}, 1 \leq i \leq r$} of $G$. Hence, we have
$$A_{b}\begin{pmatrix}
\frac{e_{\Q}(\lambda)}{p^r} & \cdots & \frac{e_{\Q}(\lambda)}{p^r}\\
\vdots & \ddots & \vdots \\
\frac{e_{\Q}(\lambda)}{p^r} & \cdots & \frac{e_{\Q}(\lambda)}{p^r}
\end{pmatrix}	=\begin{pmatrix}
e_{\Q}(\lambda) & 0 & \cdots & 0 \\
0 & 0 & \cdots & 0\\
\vdots & \vdots & \ddots & \vdots \\
0 & 0 & \cdots & 0
\end{pmatrix} A_{b}.
$$ Consequently,
\begin{enumerate}[(a)]
	\item $a_{11}=a_{12}=\cdots = a_{1p^r}=\alpha$ (say).
	\item $a_{i1}+a_{i2}+\cdots +a_{ip^r}=0$, for all $2 \leq i \leq p^r.$
\end{enumerate}	Next, we use equation (\ref{eq12}) for $k=2.$ For that, we need the expression of $b^{-1}B_{2}b$. Let $b^{-1}B_{2}b={X_{2}}$ and let $d_{ij}$ be its $(i,j)$-th entry. Again, from the expression of $\theta$ and the group relations, we can see that $d_{ij}=\frac{x^{-(j-1)}e_{\Q}(\lambda)}{p^r}$. Therefore, $${X_{2}}=\begin{pmatrix}
\frac{e_{\Q}(\lambda)}{p^r} & \frac{x^{-1}e_{\Q}(\lambda)}{p^r} & \cdots & \frac{xe_{\Q}(\lambda)}{p^r}\\
\frac{e_{\Q}(\lambda)}{p^r} & \frac{x^{-1}e_{\Q}(\lambda)}{p^r} & \cdots& \frac{xe_{\Q}(\lambda)}{p^r} \\
\vdots & \vdots & \ddots&\vdots \\
\frac{e_{\Q}(\lambda)}{p^r} & \frac{x^{-1}e_{\Q}(\lambda)}{p^r} & \cdots& \frac{xe_{\Q}(\lambda)}{p^r}
\end{pmatrix}.$$ All the rows of $X_{2}$ are {the} same and in each row of $X_2$, each power of $x$ appears precisely $p^{r-1}$ times. Consequently, from equation (\ref{eq12}),
$$A_{b}\begin{pmatrix}
\frac{e_{\Q}(\lambda)}{p^r} & \frac{x^{-1}e_{\Q}(\lambda)}{p^r} & \cdots & \frac{xe_{\Q}(\lambda)}{p^r}\\
\frac{e_{\Q}(\lambda)}{p^r} & \frac{x^{-1}e_{\Q}(\lambda)}{p^r} & \cdots& \frac{xe_{\Q}(\lambda)}{p^r} \\
\vdots & \vdots & \ddots&\vdots \\
\frac{e_{\Q}(\lambda)}{p^r} & \frac{x^{-1}e_{\Q}(\lambda)}{p^r} & \cdots& \frac{xe_{\Q}(\lambda)}{p^r}
\end{pmatrix}=\begin{pmatrix}
0 & e_{\Q}(\lambda) & 0 & \cdots & 0\\
0 & 0 & 0 & \cdots & 0\\
\vdots & \vdots & \vdots & \ddots & \vdots \\
0 & 0 & 0 & \cdots & 0
\end{pmatrix}A_{b}.$$ In view of (a), it thus turns out that the second row of $A_{b}$ is $\alpha p^r$ times the first row of $b^{-1}B_{2}b$, i.e., $a_{2j}=\alpha x^{-(j-1)}e_{\Q}(\lambda)$.\para \noindent Using equation (\ref{eq12}) for $B_{3},B_{4},\ldots,$ we can find all the rows of $A_{b}$ and see that $$A_{b}=\begin{pmatrix}
\alpha e_{\Q}(\lambda) & \alpha e_{\Q}(\lambda) & \cdots & \alpha e_{\Q}(\lambda)\\
\alpha e_{\Q}(\lambda) & {\alpha} x^{-1}e_{\Q}(\lambda) & \cdots  & {\alpha} xe_{\Q}(\lambda) \\
\vdots & \vdots & \ddots & \vdots\\
\alpha e_{\Q}(\lambda) & {\alpha} xe_{\Q}(\lambda)  & \cdots & {\alpha} x^{-r}e_{\Q}(\lambda) 
\end{pmatrix}.$$ In all rows of $A_{b}$, except first row, all powers of $x$ occur precisely $p^{r-1}$ times. By replacing $A_{b}$ by $A_{b}\operatorname{diag}(\alpha^{-1},\cdots, \alpha^{-1})_{p^r}$, we can assume that $\alpha=1.$ This is because $\operatorname{diag}(\alpha^{-1},
\ldots, \alpha^{-1}) \in \operatorname{Cen}_{\mathcal{A}}(\mathcal{B}).$ Hence, the following $A_{b}$ serves our purpose, i.e., $bA_{b} \in \operatorname{Cen}_{\mathcal{A}}(\mathcal{B})$: \begin{equation}{\label{eq13}}
A_{b}=\begin{pmatrix}
e_{\Q}(\lambda) & e_{\Q}(\lambda) & \cdots & e_{\Q}(\lambda)\\
e_{\Q}(\lambda) & x^{-1}e_{\Q}(\lambda) & \cdots &  xe_{\Q}(\lambda) \\
\vdots & \vdots & \ddots & \vdots\\
e_{\Q}(\lambda) &  xe_{\Q}(\lambda) & \cdots&  x^{-r}e_{\Q}(\lambda) 
\end{pmatrix}.
\end{equation} Next, we need to compute $z_{b}^2=bA_{b}bA_{b}.$ Note that the conjugation on $\mathbb{E}$ by $bA_{b}$ and $b(b^{-1}A_{b}b)=A_{b}b$ are both automorphisms of $\mathbb{E}/\mathbb{F}$ of order $2$. Hence, $A_{b}^{-1}b^{-1}A_{b}b \in \mathbb{E},$ i.e., \begin{equation}{\label{eq14}}
b^{-1}A_{b}b=A_{b}\operatorname{diag}(\beta, \ldots , \beta)_{p^r},~\mbox{for some}~\beta \in \Q He_{\Q}(\lambda).
\end{equation} Using the map given in isomorphism $\theta$, equation (b) together with the group relations and the fact that $(1+x+\cdots+x^{p-1})e_{\Q}(\lambda)=0$, we see that the $(1,1)$ entry of $b^{-1}A_{b}b$ is $e_{\Q}(\lambda)$. Hence, comparing the $(1,1)$ entry in equation (\ref{eq14}), we find that $\beta =e_{\Q}(\lambda),$ i.e., $b^{-1}A_{b}b=A_{b}.$ Consequently, $z_{b}^2=bA_{b}bA_{b}=b^{-1}A_{b}bA_{b}=A_{b}^2.$ From equation (\ref{eq13}), note that $(1,1)$ entry of $z_{b}^2$ is clearly $p^re_{\Q}(\lambda).$ But $z_{b}^2$ being in $\mathbb{F},$ we must have \begin{equation}{\label{eq15}}
z_{b}^2=\operatorname{diag}(p^re_{\Q}(\lambda),\ldots,p^re_{\Q}(\lambda))_{p^r}.
\end{equation} We also have $\mathbb{E} \cong \Q He_{\Q}(\lambda)$ and $\Q He_{\Q}(\lambda)$ is isomorphic to $\Q(\zeta_{2^jp})$ via the map $xe_{\Q}(\lambda) \mapsto \zeta_{p},~y_{i}e_{\Q}(\lambda)\mapsto 1$ and $ae_{\Q}(\lambda) \mapsto \zeta_{2^j}.$ Also $\mathbb{F}=\mathcal{Z}(\Q Ge_{\Q}(\lambda^G))$ is isomorphic to the subfield $F=\Q(\zeta_{2^jp}+\zeta_{2^jp}^{-1})$ of $\Q(\zeta_{2^jp})$ kept fixed by $\langle \sigma \rangle,$ where $\sigma: \Q(\zeta_{2^jp}) \rightarrow \Q(\zeta_{2^jp})$ sends $\zeta_{2^jp}$ to $\zeta_{2^jp}^{-1}.$ Consequently, equation (\ref{eq15}) together with equation (\ref{nn}) implies that $$\Q Ge_{\Q}(\lambda^G) \cong M_{p^r}(\Q(\zeta_{2^jp})/\Q(\zeta_{2^jp}+\zeta_{2^jp}^{-1}),\langle \sigma \rangle, p^r).$$ As $p \equiv 1 \operatorname{mod} 4,$ we have by (\cite{Kha}, Chapter 2, Exercise 15) that $\sqrt{p} \in F$ and hence $p \in N_{\Q(\zeta_{2^jp})/F}(\Q(\zeta_{2^jp})^{\times})$, where $\Q(\zeta_{2^jp})^{\times}$ are non-zero elements of $\Q(\zeta_{2^jp}).$ Consequently, $p^r \in N_{\Q(\zeta_{2^jp})/F}(\Q(\zeta_{2^jp})^{\times})$. This gives $$(\Q(\zeta_{2^jp})/\Q(\zeta_{2^jp}+\zeta_{2^jp}^{-1}),\langle \sigma \rangle, p^r) \cong (\Q(\zeta_{2^jp})/\Q(\zeta_{2^jp}+\zeta_{2^jp}^{-1}),\langle \sigma \rangle, 1)$$ and hence $\Q Ge_{\Q}(\lambda^G) \cong M_{2p^r}(\Q(\zeta_{2^jp}+\zeta_{2^jp}^{-1}))$. This completes the proof. ~\qed

\para \noindent \textbf{Proof of Theorem~\ref{thm8}.} From Lemmas \ref{l4} - \ref{l8}, we have that  
\begin{flalign*}
	 &~ \Q^{(4)} \bigoplus_{i=2}^{n-1} M_{2}(\Q(\zeta_{2^i}+\zeta_{2^i}^{-1})) \bigoplus M_{2^{n-1}}(\Q(\zeta_{p}+\zeta_{p}^{-1}))^{(4m_{1})}\\ & \bigoplus  M_{2^n}(\Q(\zeta_{p}+\zeta_{p}^k+\cdots + \zeta_{p}^{k^{2^{n-1}-1}}))^{(m_{1})}\bigoplus M_{2^{n}}(\Q(\zeta_{p}+\zeta_{p}^{-1}))^{(m_{2})}\\ & \bigoplus_{i=0}^{n-1} M_{2p^r}(\Q(\zeta_{2^ip}+\zeta_{2^ip}^{-1})) ,
	\end{flalign*}
 {is a subalgebra of $\Q G$, where  $m_{1}=\frac{p^r-1}{p-1}$ and $m_{2}=\frac{(p^{2r}-1)-(2^{n-1}+2)(p^r-1)}{(p-1)2^{n-1}}$. If $\zeta$ is $m^{th}$ root of unity and $\phi$ is the Euler phi function, then we know that $\operatorname{dim}_{\mathbb{Q}}(\Q(\zeta+\zeta^{-1}))=\frac{\phi(m)}{2}.$ Also, from (\cite{DF04}, Section 14.5, example (2)), we have $\operatorname{dim}_{\Q}(\Q(\zeta_{p}+\zeta_{p}^k+\cdots + \zeta_{p}^{k^{2^{n-1}-1}}))=\frac{p-1}{2^{n-1}}=t,$ because $p-1=2^{n-1}t$. Hence, the dimension over $\Q$ of the above subalgebra is equal to $$4+4\sum_{i=2}^{n-1}\frac{\phi(2^i)}{2}+4m_{1}2^{2n-2}\frac{\phi(p)}{2}+m_{1}2^{2n}t+m_{2}2^{2n}\frac{\phi(p)}{2}+4p^{2r}\sum_{i=0}^{n-1}\frac{\phi(2^ip)}{2}=p^{2r+1}2^n.$$As the order of $|G|$ is $p^{2r+1}2^n$, we obtain that the subalgebra given above is equal to $\Q G$. }~\qed 
\subsection{Primitive Idempotents}
In this section, we illustrate Theorem~\ref{thm3} by providing explicit computation of the orthogonal primitive idempotents in each simple component of $\Q G$. We follow the cases from  Lemmas \ref{l4} - \ref{l8} to list the simple components given by (generalized) strong Shoda pairs $(H,K)$.
\subsubsection{\bf $(H,K)\in \left\{ (G,G), (G,\langle P,a^2,ab \rangle), (G,\langle P,a \rangle), (G,\langle P,a^2,b \rangle)\right\}$.} These are the possibilities for $(H,K)$ with $H=G$ from Section 5.2. In all these cases, $\Q G e(G,H,K)$ is isomorphic to $\Q$ and $\varepsilon(H,K)$ is the corresponding primitive idempotent.
\subsubsection{$(H,K)$ with $H=\langle P,a \rangle$ and $K = \langle P,a^{2^i} \rangle,~2 \leq i \leq n-1$.}
These are strong Shoda pairs {from Lemma \ref{l4}.} Recall that in this case $N_{G}(K)=G$ and $$\Q Ge(G,H,K)\cong \Q H\varepsilon(H,K)*^{\sigma}G/H,$$ with $\sigma$ given by $bH \mapsto \sigma_{b},$ where $\sigma_{b}(a\varepsilon)=b^{-1}a\varepsilon b=a^{-1}\varepsilon$ and $\varepsilon=\varepsilon(H,K).$ Consequently, $\mathbb{F}=F=\Q(a\varepsilon+a^{-1}\varepsilon),~\mathcal{B}=\mathbb{F},~\operatorname{Cen}_{\mathcal{A}}(\mathcal{B})=\mathcal{A}$ and $\mathbb{E}=E=\Q H\varepsilon$.\para \noindent In view of Remark~(\ref{rr},~a), $\{1,b\}$ is a {basis of units} for $\operatorname{Cen}_{\mathcal{A}}(\mathcal{B})$ (as a vector space over $E$) and we can take $\mathbf{z}_{1}=1,~\mathbf{z}_{2}=b$ and $\widehat{E}=\frac{1+b}{2}$.\para \noindent As $[E:F]=2$ and $\{(1+a)\varepsilon,\sigma_{b}((1+a)\varepsilon)\}$ are linearly independent over $F,$ we have $w=(1+a)\varepsilon$ is a normal element of extension $E/F.$\para \noindent Another ingredient required to apply Theorem~\ref{thm3} is $\alpha= \alpha_{0}+\alpha_{1}b,$ where $\alpha_{0},\alpha_{1} \in \Q H\varepsilon$ satisfy the following system of equations:
\begin{flalign*}
\begin{pmatrix}
(1+a)\varepsilon & (1+a^{-1})\varepsilon \\
(1+a^{-1})\varepsilon & (1+a)\varepsilon
\end{pmatrix} \begin{pmatrix}
\alpha_{0}\\
\alpha_{1}
\end{pmatrix}=\begin{pmatrix}
a\varepsilon+a^{-1}\varepsilon +2\varepsilon\\
a\varepsilon-a^{-1}\varepsilon
\end{pmatrix}.
\end{flalign*}
Solving the above system of equations, we get $\alpha=\frac{a^2\varepsilon+\varepsilon}{a^2\varepsilon-\varepsilon}+ \frac{a^2\varepsilon-2a\varepsilon-\varepsilon}{a^2\varepsilon-\varepsilon} b.$ \para \noindent All of these ingredients yield that
$$\{ \alpha^{-1}\widehat{E}\varepsilon \alpha, b^{-1}\alpha^{-1}\widehat{E}\varepsilon \alpha b \}$$ are orthogonal primitive idempotents of $\Q G$ corresponding to the simple component $\Q Ge(G,H,K)$.
\subsubsection{$(H,K)$ with $H=\langle P,b \rangle$ and $K \lhd H$, $H/K$ cyclic, $K \neq P$.} These are strong Shoda pairs from Lemma \ref{l5} for $H=\langle P,b \rangle$. If $hK$ is a generator of $H/K,$ we have $$\Q Ge(G,H,K)\cong M_{2^{n-2}}(\Q H\varepsilon(H,K)*^{\sigma} \langle P,b,a^{2^{n-2}}\rangle/H),$$ with $\sigma$ given by $a^{2^{n-2}}H \mapsto\sigma_{a^{2^{n-2}}}$, where $\sigma_{a^{2^{n-2}}}(h\varepsilon)=a^{-2^{n-2}}h\varepsilon a^{2^{n-2}}$ and $\varepsilon(H,K)=\varepsilon.$ We now have $F=\Q(h\varepsilon+a^{-2^{n-2}}ha^{2^{n-2}}\varepsilon),~\mathcal{B}\cong M_{2^{n-2}}(F),~\operatorname{Cen}_{\mathcal{A}}(\mathcal{B}) \cong \Q H \varepsilon*^{\sigma}\langle a^{2^{n-2}}H \rangle,~E=\Q H \varepsilon,~T=\{1,a,\ldots,a^{2^{n-2}-1}\},$ $\mathbf{z}_{1}=1,~\mathbf{z}_{2}=a^{2^{n-2}}$ and $\widehat{E}=\frac{1+a^{2^{n-2}}}{2}.$ Furthermore, as $K\neq P,$ we can choose $i$ such that $y_{i}\notin K.$ As $[E:F]=2$ and $\{y_{i}\varepsilon,\sigma_{a^{2^{n-2}}}(y_{i}\varepsilon)\}$ are linearly independent over $F,$ we have $w=y_{i}\varepsilon$ is a normal element of extension $E/F$ and Theorem~\ref{thm3} thus yield that \begin{flalign*}
\{ a^{-k}\alpha^{-1}\widehat{E}\varepsilon\alpha a^{k},a^{-k-2^{n-2}}\alpha^{-1}\widehat{E}\varepsilon \alpha a^{k+2^{n-2}}~|~0\leq k\leq 2^{n-2}-1 \}
\end{flalign*} are orthogonal primitive idempotents of $\Q Ge(G,H,K)$, where $\alpha= \alpha_{0}+\alpha_{1}a^{2^{n-2}}$ with $\alpha_{0},\alpha_{1} \in \Q H\varepsilon$ satisfying the following system of equations: \begin{flalign*}
\begin{pmatrix}
y_{i}\varepsilon & a^{-2^{n-2}}y_{i}a^{2^{n-2}}\varepsilon \\
a^{-2^{n-2}}y_{i}a^{2^{n-2}}\varepsilon & y_{i}\varepsilon
\end{pmatrix} \begin{pmatrix}
\alpha_{0}\\
\alpha_{1}
\end{pmatrix}=\begin{pmatrix}
y_{i}\varepsilon+a^{-2^{n-2}}y_{i}a^{2^{n-2}}\varepsilon\\
y_{i}\varepsilon-a^{-2^{n-2}}y_{i}a^{2^{n-2}}\varepsilon
\end{pmatrix}.
\end{flalign*}
\subsubsection{$(H,K)$ with $H=\langle P,ba \rangle$ and  $K\lhd H$, $H/K$ cyclic, $K \neq P$.}
These are strong Shoda pairs from Lemma~\ref{l5} for  $H=\langle P,ba \rangle$. This case is similar to the above and we have that 
\begin{flalign*}
\{ a^{-k}\alpha^{-1}\widehat{E}\varepsilon\alpha a^{k},a^{-k-2^{n-2}}\alpha^{-1}\widehat{E}\varepsilon \alpha a^{k+2^{n-2}}~|~0 \leq k \leq 2^{n-2}-1 \}
\end{flalign*} are orthogonal primitive idempotents of $\Q Ge(G,H,K)$. Here, $\widehat{E}=\frac{1+a^{2^{n-2}}}{2},~\varepsilon=\varepsilon(H,K)$ and $\alpha= \alpha_{0}+\alpha_{1}a^{2^{n-2}}$ with $\alpha_{0},\alpha_{1} \in \Q H\varepsilon$ satisfying the following system of equations: \begin{flalign*}
\begin{pmatrix}
y_{i}\varepsilon & a^{-2^{n-2}}y_{i}a^{2^{n-2}}\varepsilon \\
a^{-2^{n-2}}y_{i}a^{2^{n-2}}\varepsilon & y_{i}\varepsilon
\end{pmatrix} \begin{pmatrix}
\alpha_{0}\\
\alpha_{1}
\end{pmatrix}=\begin{pmatrix}
y_{i}\varepsilon+a^{-2^{n-2}}y_{i}a^{2^{n-2}}\varepsilon\\
y_{i}\varepsilon-a^{-2^{n-2}}y_{i}a^{2^{n-2}}\varepsilon
\end{pmatrix}.
\end{flalign*} Here $i$ is such that $y_{i}\notin K.$
\subsubsection{$(H,K)$ with $H=P$ and  $K \unlhd H$, $H/K$ cyclic, $K\in Z$.}
This is case (i) from Lemma~\ref{l6} with $Z$ defined as in equations (\ref{eqo}). 
 In this case, we obtain the following primitive idempotents for $\Q Ge(G,H,K)$: \begin{flalign*}
\{b^{-j}a^{-i}\alpha^{-1}\widehat{E}\varepsilon \alpha a^{i}b^{j}~|~0 \leq i \leq 2^{n-1}-1,~j = 0,1\},
\end{flalign*} where $\widehat{E}=\frac{1+a+\cdots+a^{2^{n-1}-1}}{2^{n-1}},~\varepsilon=\varepsilon(H,K),~\langle hK \rangle$ is a generator of $H/K$ and $\alpha= \alpha_{0}+\alpha_{1}a+\cdots+\alpha_{2^{n-1}-1} a^{2^{n-1}-1},$ where $\alpha_{i} \in \Q H\varepsilon(H,K)$ satisfy 
\begin{footnotesize}
\begin{flalign*}
\begin{pmatrix}
h\varepsilon & a^{-1}ha\varepsilon & \cdots & a^{-(2^{n-1}-1)}ha^{2^{n-1}-1}\varepsilon \\
a^{-1}ha\varepsilon & a^{-2}ha^2\varepsilon & \cdots & h\varepsilon \\
\vdots & \vdots & \ddots & \vdots \\
a^{-(2^{n-1}-1)}ha^{2^{n-1}-1}\varepsilon & h\varepsilon  & \cdots&  a^{-(2^{n-1}-2)}ha^{2^{n-1}-2}\varepsilon  
\end{pmatrix} \begin{pmatrix}
\alpha_{0}\\
\alpha_{1}\\
\vdots \\
\alpha_{2^{n-1}-1}
\end{pmatrix}=\begin{pmatrix}
\sum_{k=0}^{2^{n-1}-1}a^{-k}ha^{k}\varepsilon\\
h\varepsilon-a^{-1}ha\varepsilon\\
\vdots \\
h\varepsilon -aha^{-1}\varepsilon
\end{pmatrix}.
\end{flalign*}
\end{footnotesize}
\subsubsection{$(H,K)$ with $H=P$ and  $K \unlhd H$, $H/K$ cyclic,  $K \notin Y$ and $K \notin Z$.}  This is case (ii) from Lemma~\ref{l6} with $Y$ and $Z$ defined as in equations (\ref{eq6}) and (\ref{eqo}). The orthogonal primitive idempotents of $\Q Ge(G,H,K)$ are:
\begin{flalign*}
\{b^{-j}a^{-i}\alpha^{-1}\widehat{E}\varepsilon \alpha a^ib^{j},b^{-j}a^{-i-2^{n-2}}\alpha^{-1}\widehat{E}\varepsilon \alpha a^{2^{n-2}+i}b^{j}~|~0 \leq i \leq 2^{n-2}-1,~0 \leq j \leq 1\}.
\end{flalign*} Here, $\widehat{E}=\frac{1+a^{2^{n-2}}}{2},~\varepsilon=\varepsilon(H,K),~\langle hK\rangle$ is a generator of $H/K$ and $\alpha= \alpha_{0}+\alpha_{1}a^{2^{n-2}},$ with $\alpha_{0},\alpha_{1} \in \Q H\varepsilon(H,K)$ satisfying 
\begin{flalign*}
\begin{pmatrix}
h\varepsilon & h^{-1}\varepsilon \\
h^{-1}\varepsilon & h\varepsilon
\end{pmatrix} \begin{pmatrix}
\alpha_{0}\\
\alpha_{1}
\end{pmatrix}=\begin{pmatrix}
h\varepsilon+h^{-1}\varepsilon\\
h\varepsilon-h^{-1}\varepsilon
\end{pmatrix}.
\end{flalign*}
\subsubsection{$(H,K)$ with $H=\langle x,y_{1},\ldots,y_{r},a \rangle$ and $K=\langle y_{1},\ldots,y_{r},a^{2^j}\rangle$, for $0 \leq j \leq n-1$.}
In Lemma~\ref{l7}, we have seen that every such $(H,K)$  is a generalized strong Shoda pair with $$H=H_0 \leq H_{1}=\langle P,a \rangle \leq G$$ as a strong inductive chain and the following isomorphism holds $$\mathcal{A}\cong\Q Ge_{\Q}(\lambda^G) \cong M_{p^r}(\Q H\varepsilon(H,K)*\langle \sigma_{b}\rangle).$$ Here, we have $E=\Q H\varepsilon(H,K),~\mathbb{E}=\{\operatorname{diag}(\alpha,\ldots,\alpha)_{p^r}~|~\alpha \in E\},~\mathbb{F}=\mathcal{Z}(\Q Ge_{\Q}(\lambda^G)),$ $F=(\Q H\varepsilon(H,K))^{\langle\sigma_{b}\rangle},~\mathcal{B}=M_{p^r}(F),~\operatorname{Cen}_{\mathcal{A}}(\mathcal{B})=\Q H\varepsilon(H,K)*\langle \sigma_{b}\rangle$ and $T=T_{0}=\{z_{1}^{i_{1}}\cdots z_{r}^{i_{r}}~|~0 \leq i_{j} \leq p-1\}$, a transversal of $H$ in $H_{1}=\langle P,a \rangle$. In Lemma~\ref{l8}, we have seen that $~\mathbf{z}_{1}=1,~\mathbf{z}_{2}=z_{b}=bA_{b}$ is a {basis of units} of $\operatorname{Cen}_{\mathcal{A}}(\mathcal{B})$ as $\mathbb{E}$-vector space. Since $z_{b}^2=p^re_{\Q}(\lambda),$ {replacing $\mathbf{z}_{2}$ by $p^{-r/2}\mathbf{z}_{2}$, we can assume that $\mathbf{z}_{1}=1$ and $\mathbf{z}_{2}=p^{-r/2}bA_{b}$ as desired in equation (\ref{eq2}).} \para \noindent As $[E:F]=2$ and $\{ax\varepsilon,\sigma_{b}(ax\varepsilon)\}$ being linearly independent over $F$, it turns out that $w=ax\varepsilon$ is a normal element of $E/F.$ Consequently, we obtain the following set of orthogonal primitive idempotents of $\mathcal{A}=\Q Ge_{\Q}(\lambda^G)$:
\begin{flalign*}
\{ \gamma^{-1}\mathbf{z}_{i}^{-1}\alpha^{-1}\widehat{E}\varepsilon(H,K) \alpha \mathbf{z}_{i} \gamma~|~i=1,2 \},
\end{flalign*} where $\mathbf{z}_{1}=1,~\mathbf{z}_{2}=p^{-r/2}bA_{b}$, $\widehat{E}=\frac{\mathbf{z}_{1}+\mathbf{z}_{2}}{2}$, $\gamma \in T$ and $$\alpha=\frac{(ax)^2\varepsilon+(ax)^{-2}\varepsilon}{(ax)^2\varepsilon-(ax)^{-2}\varepsilon}\mathbf{z}_{1}+ \frac{(ax)^2\varepsilon-(ax)^{-2}\varepsilon-2\varepsilon}{(ax)^2\varepsilon-(ax)^{-2}\varepsilon} \mathbf{z}_{2}.$$
\begin{remark}  Notice that the computations for this example are quite laborious and complex, especially the ones for the (generalized) strong Shoda pairs. However, for some concrete values of the prime $p$, one can use GAP~\cite{GAP} and the GAP package Wedderga~\cite{Wedderga} in order to compute at least the strong Shoda pairs of $G$ and use that knowledge to compute the primitive idempotents in the corresponding simple components using the theory developed. These GAP computations can be done only for relative small values of the prime $p$ that are computationally supported by the GAP system.
\end{remark}
\subsection{Units}
In the Wedderburn decomposition of $\Q G$, there is one exceptional component, namely $M_{2}(\Q)$, arising from the Shoda pair $(\langle P,a\rangle,\langle P,a^4 \rangle)$ (obtained from Lemma \ref{l4} with $i=2$) and  all other components are non exceptional. Notice that this exceptional component is coming from the Wedderburn decomposition of $\Q (G/P)\cong \Q \langle a,b\rangle\cong \Q D_{2^n}$ and Ritter and Sehgal (\cite{Seh93}, Theorem 23.1) proved that the Bass cyclic and the bicyclic units of $\mathbb{Z}D_{2^n}$ generate a subgroup of finite index in $\mathcal{U}(\mathbb{Z}D_{2^n})$. In this section, we will use Ritter and Sehgal's idea of working with exceptional component of $\Q D_{2^n}$ together with Theorem \ref{thm7} (for non exceptional components of $\Q G$) to produce a subgroup of finite index in $\mathcal{U}(\mathbb{Z}G).$ \para \noindent Let $X$ be the set of all the Shoda pairs of $G$ obtained in Section 5.2 except $(H,K)=(\langle P,a\rangle,\langle P,a^4 \rangle).$ Let $E_X$ be the set of primitive central idempotents of $\Q G$ corresponding to the Shoda pairs in $X$.\para \noindent The following table describes the subgroup $\langle V_{e}^{+},V_{e}^{-}\rangle$ of $\mathcal{U}(\mathbb{Z}G)$ for all $e \in E_X.$
\para\noindent 
\begin{scriptsize}\begin{tabular}{|c|c|c|c|}
	\hline  \rule[-2ex]{0pt}{4.5ex}  $(H,K)$  & Matrix units of  &$\mathtt{B}$ &$\langle V_{e}^{+},V_{e}^{-}\rangle$ \\ 
	 & $\Q Ge$ & & \\ \hline
	 	\hline  \rule[-2ex]{0pt}{4.5ex}  $(G,K),$  & $\varepsilon$ & $\{\varepsilon\}$&   identity subgroup \\
	 	\rule[-2ex]{0pt}{4.5ex}  $K\in \{G,\langle P,a^2,ab\rangle$,  &  & &   \\  $\langle P,a^2,b\rangle,\langle P,a\rangle\}$	 &  & & \\ \hline  \rule[-2ex]{0pt}{4.5ex}  $(\langle P,a\rangle,\langle P,a^{2^i}\rangle),$  &$b^{-i}\alpha^{-1}\widehat{E}_{1}\alpha b^j,$ & $\{(a^l+a^{-l}+2)\varepsilon|$&  $\langle 1+c\beta b^{-i}\alpha^{-1}\widehat{E_{1}}\alpha b^j|$ \\
	 	  \rule[-2ex]{0pt}{4.5ex} $3\leq i\leq n$   & $0\leq i,j\leq 1$ & $l$ odd, $1 \leq l\leq 2^{i-2}\}$ &$\beta \in \mathtt{B},0\leq i\neq j\leq 1\rangle$   \\
	\hline  \rule[-2ex]{0pt}{4.5ex}  $(\langle P,b \rangle,K),$  & $a^{-i}\alpha^{-1}\widehat{E_{2}}\varepsilon \alpha a^j$ & $\{y_{i}^{l}\varepsilon+y_{i}^{-l}\varepsilon~|$&  $\langle 1+c \beta a^{-i}\alpha^{-1}\widehat{E_{2}}\varepsilon \alpha a^{j} |$  \\$K\in X\backslash \{P\}$&   $0 \leq i,j \leq 2^{n-1}-1$& $1\leq l \leq \frac{p-1}{2}\},$  & $ \beta \in \mathtt{B}, 0 \leq i\neq j \leq 2^{n-1}-1\rangle$  \\& & where $i$ is such that $y_{i}\notin K$. &  \\
	\hline  \rule[-2ex]{0pt}{4.5ex}  $(\langle P,ba \rangle,K),$  & $a^{-i}\alpha^{-1}\widehat{E_{2}}\varepsilon \alpha a^j$ & $\{y_{i}^{l}\varepsilon+y_{i}^{-l}\varepsilon~|$& $\langle 1+c \beta a^{-i}\alpha^{-1}\widehat{E_2}\varepsilon \alpha a^{j} |$   \\$K\in X'\backslash \{P\}$&   $0 \leq i,j \leq 2^{n-1}-1$&  $1\leq l \leq \frac{p-1}{2}\},$ & $ \beta \in \mathtt{B}, 0 \leq i\neq j \leq 2^{n-1}-1\rangle$ \\ & & where $i$ is such that $y_{i}\notin K$. &  \\
	\hline \rule[-2ex]{0pt}{4.5ex}  $(P,K),~K \in Z$ & $	b^{-i_{1}}a^{-i}\alpha^{-1}\widehat{E_{3}}\varepsilon \alpha a^{j}b^{j_{1}},$ & $\{\delta_{y_u}^{(m)}\varepsilon~|~m\in \mathcal{T}_{1}\},$ &  $\langle 1+c \beta b^{-i_{1}}a^{-i}\alpha^{-1}\widehat{E_{3}}\varepsilon \alpha a^{j}b^{j_{1}} | $   \\ 
	  \rule[-2ex]{0pt}{4.5ex} &   $0 \leq i, j \leq 2^{n-1}-1$ &  if  $P/K=\langle y_{u}K\rangle$ or &  $ \beta \in \mathtt{B}, 0 \leq i\neq j \leq 2^{n-1}-1,$ \\
	&  $0 \leq i_{1}, j_{1} \leq 1$ &   $\{\delta_{z_u}^{(m)}\varepsilon~|~m \in \mathcal{T}_{2}\},$   & $0 \leq i_{1}\neq j_{1} \leq 1\rangle$ \\
	 \rule[-2ex]{0pt}{4.5ex} & &    if $P/K=\langle z_{u}K\rangle$  & \\
	\hline \rule[-2ex]{0pt}{4.5ex}  $(P,K),~ K \notin Y,$ & $b^{-i_{1}}a^{-i}\alpha^{-1}\widehat{E_{2}}\varepsilon \alpha a^{j}b^{j_{1}}$, &$\{y_{i}^{l}\varepsilon+y_{i}^{-l}\varepsilon~|$  & $\langle 1+c \beta b^{-i_{1}}a^{-i}\alpha^{-1}\widehat{E_{2}}\varepsilon \alpha a^{j}b^{j_{1}} | $  \\ 
	 \rule[-2ex]{0pt}{4.5ex} $K \notin Z$&  $0 \leq i, j \leq 2^{n-1}-1$& $1\leq l \leq \frac{p-1}{2}\},$ &   $ \beta \in \mathtt{B}, 0 \leq i\neq j \leq 2^{n-1}-1,$ \\&   $0 \leq i_{1}, j_{1} \leq 1$& where $i$ is such that $y_{i}\notin K$. &$0 \leq i_{1}\neq j_{1} \leq 1\rangle$   \\
	\hline  \rule[-2ex]{0pt}{4.5ex} $(\langle x,y_{1},\ldots,y_{r},a \rangle,$ & $\gamma^{-1}\mathbf{z}_{i}^{-1}\alpha^{-1}\widehat{E_4}\varepsilon \alpha \mathbf{z}_{i'} \gamma'$,  &$\{(ax)^{l}\varepsilon+(ax)^{-l}\varepsilon~|$ & $\langle 1+c \beta \gamma^{-1}\mathbf{z}_{i}^{-1}\alpha^{-1}\widehat{E_4}\varepsilon \alpha \mathbf{z}_{i'} \gamma' ~|$  \\ $\langle y_{1},\cdots,y_{r},a^{2^j}\rangle)$& $1\leq i,i' \leq 2,~\gamma,\gamma' \in T$ & $~1\leq l\leq 2^{j-2}(p-1)\}$& $ \beta \in \mathtt{B},~1\leq i\neq i' \leq 2\rangle$\\$0 \leq j \leq n-1$ & & & \\
		\hline 
\end{tabular}
\end{scriptsize}
\para\noindent
We use the notation $\varepsilon=\varepsilon(H,K),~\delta_{y_u}^{(m)}= y_{u}^{m}+y_{u}^{km}+\cdots+y_{u}^{k^{2^{n-1}-1}m}$, $\delta_{z_u}^{(m)}=z_{u}^{m}+z_{u}^{qm}+\cdots+z_{u}^{q^{2^{n-1}-1}m},$ $\mathcal{T}_{1},~\mathcal{T}_{2}$
transversals of $\langle k\rangle,\langle q\rangle$ in $\mathbb{Z}_{p}^{\times}$
respectively, where $\mathbb{Z}_{p}^{\times}$
is a unit group of $\mathbb{Z}/{p\mathbb{Z}}$, $\widehat{E_1}=\frac{1+b}{2}$, $\widehat{E_{2}}=\frac{1+a^{2^{n-2}}}{2}$, $\widehat{E_{3}}=\frac{1+a+\cdots+a^{2^{n-1}-1}}{2^{n-1}}$, $\widehat{E_4}=\frac{\mathbf{z}_1+\mathbf{z}_2}{2}$ and $\alpha'$s are as mentioned in the respective subsections of Section 5.3.  \para \noindent Recall that a bicyclic unit in an integral group ring $\mathbb{Z}G$ is of unit of the type $u_{g,\tilde{h}}=1+(1-h)g\tilde{h}$, where $g$ and $h$ are elements of finite group $G$ and $\tilde{h}=1+h+\cdots+h^{|h|-1}.$ 
\begin{theorem}
Let $V=\{u_{a,\tilde{a^2b}},u_{a,\tilde{b}},u_{b,\tilde{ab}},u_{a,\tilde{a^3b}}\}$ and let $U$ be a subgroup of finite index in $\mathcal{Z}(\mathcal{U}(\mathbb{Z}G)).$ Then $$\mathcal{V}=\langle U,V,\prod_{e\in E_X}V_e^+,\prod_{e\in E_X}V_e^-\rangle$$ is a subgroup of finite index in $\mathcal{U}(\mathbb{Z}G).$
\end{theorem}
\noindent \textbf{Proof.} We have that $\Q G\cong \bigoplus_{e\in E_X}\Q Ge\bigoplus M_2(\Q).$\para \noindent \textbf{(i)} For every $e \in E_X$, the simple component $\Q Ge$ is non exceptional. From Theorem \ref{thm7}, we have that $\mathcal{V}$ contains a subgroup, namely $\langle V_e^+,V_e^-\rangle$, which is of finite index in $1-e+SL_{[G:H]}(O_e)$ for some order $O_e$ in the center of $\Q Ge$.\para \noindent \textbf{(ii)} Let $\pi$ be the projection of $\Q G$ on $M_2(\Q)$. We will prove that $\pi(\mathcal{V})$ contains a subgroup of finite index in $SL_2(\mathbb{Z}).$ The simple component $M_2(\Q)$ of $\Q G$ comes from the Wedderburn decomposition of $\Q(G/P)\cong \Q \langle a,b\rangle.$ If $T$ denotes the projection of $\Q \langle a,b\rangle$ on $M_2(\Q)$ and $S$ is the natural projection of $\Q G$ on $\Q \langle a,b\rangle$, then $\pi=T\circ S.$  Since $S$ is identity on $\Q \langle a,b\rangle,$ we have $\pi(V)=T(V).$ From part (ii) of the proof of Theorem 23.1 of \cite{Seh93}, we have that if $X_1=  T(u_{a,\tilde{a^2b}}),~X_2=T(u_{a,\tilde{b}}),~X_3=T(u_{b,\tilde{ab}}),~X_4=T(u_{a,\tilde{a^3b}})$, then $\langle X_1,X_2,X_3,X_4\rangle$ contains the subgroup $\langle X_1,X_2,X_3^2,X_4^2,X_3X_4^{-1}\rangle$, which is of finite index in $SL_2(\mathbb{Z}).$ Thus $\pi(V)$ and consequently $\pi(\mathcal{V})$ contains a subgroup of finite index in $SL_2(\mathbb{Z}).$
\para \noindent \textbf{(iii)} We now apply (\cite{Seh93}, Lemma 22.10) with $A=\bigoplus_{e \in E_X} SL_{[G:H]}(O_e)$ and $C=SL_{2}(\mathbb{Z})$. Its hypothesis is satisfied in view of $\textbf{(i)}$ and $\textbf{(ii)}$ and we obtain that $\mathcal{V}$ contains a subgroup which is of finite index in $\bigoplus_{e \in E_X} SL_{[G:H]}(O_e)\bigoplus SL_{2}(\mathbb{Z})$. Since $\mathcal{V}$ contains $U$, which is of finite index in $\mathcal{Z}(\mathcal{U}(\mathbb{Z}G))$, it turns out that the index of $\mathcal{V}$ in $\bigoplus_{e \in E_X} GL_{[G:H]}(O_e)\bigoplus GL_{2}(\mathbb{Z})$ is finite. As $\mathbb{Z}G$ and $\bigoplus_{e \in E_X} M_{[G:H]}(O_e)\bigoplus M_{2}(\mathbb{Z})$ are both orders in $\Q G$, their unit groups are commensurable and hence the index of $\mathcal{V}$ in $\mathcal{U}(\mathbb{Z}G)$ is finite.~\qed
{\section{Example 2: A class of Frobenius groups} 
\begin{theorem} If $G$ is a Frobenius group of odd order with cyclic complement, then $G$ is generalized strongly monomial, all simple components of $\Q G$ have Schur index 1 and none of them is exceptional.
\end{theorem}
\noindent\textbf{Proof.} Let $G$ be a Frobenius group with kernel $N$ and complement $K$. From (\cite{Hup}, Theorems 16.7 a) and 18.7), it is known that $N$ is nilpotent and every irreducible character of $G$ is one of the following:
\begin{enumerate}[(i)]
\item $\chi$ is lifted from an irreducible character of $K\cong G/N$.
\item $\chi$ is induced from an irreducible character of $N$.
\end{enumerate} If $\chi$ is of type (i) and $K$ is cyclic, then $\chi$ is a linear and the simple component $\Q Ge_{\Q}(\chi)$ is a field and hence, in this case, it is never exceptional and also $m_{\Q}(\chi)=1.$
\para \noindent Suppose $\chi$ is of type (ii) and it is induced from an irreducible character, say $\psi$, of $N$. As $N$ is nilpotent, there exists a strong Shoda pair $(H,K)$ of $N$ and a linear character $\lambda$ of $H$ with kernel $K$ such that $\psi=\lambda^N.$ Now, $\chi=\psi^G=(\lambda^N)^G=\lambda^G,$ i.e., $(H,K)$ is a Shoda pair of $G$ and $\chi$ arises from $(H,K)$. As $N\unlhd G,$ in view of Remark (\ref{remark3}, (i)), $(H,K)$ turns out to be a generalized strong Shoda pair of $G$ with strong inductive chain $H_{0}=H\leq H_{1}=N\leq H_{2}=G$ from $H$ to $G$.\para \noindent From (\cite{Iss}, Problem 10.2), $m_{\Q}(\chi)||G/N|m_{\Q}(\psi)$, i.e., $m_{\Q}(\chi)$ divides $|G/N|m_{\Q}(\psi)$. Also from (\cite{Iss}, Lemma 10.8), we have that $m_{\Q}(\chi)|\psi(1)$. Therefore,\linebreak $m_{\Q}(\chi)|\operatorname{gcd}(|G/N|m_{\Q}(\psi),\psi(1)).$ Since $m_{\Q}(\psi)|\psi(1)$, $\psi(1)||N|$ and $\operatorname{gcd}(|N|,|G/N|)=1$, we get $m_{\Q}(\chi)|m_{\Q}(\psi).$ Also from (\cite{Iss}, Problem 10.1), $m_{\Q}(\psi)|m_{\Q}(\chi)$. Consequently, $m_{\Q}(\chi)=m_{\Q}(\psi).$ Now, $N$ being nilpotent of odd order, by (\cite{Iss}, Corollary 10.4), $m_{\Q}(\psi)=1$ and hence $m_{\Q}(\chi)=1.$\para \noindent Let us now show that $\Q Ge_{\Q}(\chi)$ is not exceptional. From Theorem \ref{thm1}, \begin{equation}\label{neweq14}
\Q Ge_{\Q}(\chi)\cong M_{k_1}(M_{k_0}(\Q H\varepsilon*_{\tau_{0}}^{\sigma_{0}}C_{0}/H_{0})*_{\tau_{1}}^{\sigma_{1}}C_{1}/H_{1}),
\end{equation} where $k_{i}$'s, $\sigma_{i}$'s, $\tau_{i}$'s and $C_{i}$'s are as in Notation \ref{not2}.\para \noindent As $m_{\Q}(\chi)=1,$ we have \begin{equation}\label{neweq13}
\Q Ge_{\Q}(\chi)\cong M_{r}(F),
\end{equation} for some integer $r \geq 1,$ and $F$ as in Notation \ref{not2}. From (\cite{BK22}, Proposition 1), $\operatorname{dim}_{\Q}(F)=\operatorname{dim}_{\Q}(\mathcal{Z}(\Q Ge_{\Q}(\chi)))=\frac{\phi([H:K])}{|C_0/H_0||C_1/H_1|}.$\para\noindent Comparing the dimension over $\Q$ in (\ref{neweq14}) and (\ref{neweq13}), we get $r=k_0k_1|C_0/H_0||C_1/H_1|$ and hence $$\Q Ge_{\Q}(\chi)\cong M_{k_0k_1|C_0/H_0||C_1/H_1|}(F).$$If the order of $G$ is odd, then $k_0k_1|C_0/H_0||C_1/H_1|$ being divisor of $|G|$ is never 2 and hence $\Q Ge_{\Q}(\chi)$ is never exceptional.\qed \para \noindent \begin{remark}
(i) There do exist Frobenius groups of odd order with cyclic complement. The Frobenius group constructed by N. It\^o (see \cite{Hup67}, p.499, 8.6 Beispiele c)), whose Fitting subgroup has arbitrary large nilpotency class is an example of this type.\\
(ii) The smallest non-abelian Frobenius group of odd order is $C_{7}\rtimes C_3$ and it has cyclic complement, where $C_n$ is cyclic group of order $n$. Moreover $C_{q^m}\rtimes C_{p^n}$ with $p,~q$ odd primes and $C_{p^n}$ acting faithfully on $C_{q^m}$ is also a Frobenius group of odd order with cyclic complement and hence our method extends the previous one considered in \cite{JOdRVG2013}.\\
(iii) For primitive idempotents and units, the computations similar to those in Sections 5.3 and 5.4 can be done once the structure of $N$ and the action of $G$ on $N$ is known.
\end{remark}}
\section{Conclusion}
In this article, we have given a method to explicitly compute a complete set of orthogonal primitive idempotents in a simple component with Schur index 1 of a rational group algebra $\mathbb{Q}G$ for $G$ a finite generalized strongly monomial group.

The main motivation for the computation of such a complete set of primitive idempotents was the application of this information to construct units of the integral group ring $\mathbb{Z}G$ for $G$ as above and  with no exceptional simple components in  $\mathbb{Q}G$. 

The computation of a complete set of orthogonal primitive idempotents in a simple component of a  group algebra can also be of interest for a semisimple finite group algebra $\mathbb{F}G$ and $G$ a finite generalized strongly monomial group. The knowledge of such a set of primitive idempotents can be used afterwards to compute non-abelian left group codes for this class of groups $G$. We refer to \cite{OVG2015} for the development in this direction.

\section*{Acknowledgments}

The authors would like to thank the referees for the many constructive comments and remarks that helped us to improve the article.

\end{document}